\documentclass[11pt]{article}

\input{amssym.def}
\input{amssym}
\usepackage{eufrak}

\newtheorem{theorem}{Theorem}

\newtheorem{definition}[theorem]{Definition}

\newtheorem{proposition}[theorem]{Proposition}
\newtheorem{remark}[theorem]{Remark}

\def\K{{\Bbb K}}

\def\ot{{\otimes}}

\def\End{{\rm End}}
\def\vv{V^{\otimes 2}}

\def\sscr#1{{\!\raisebox{-1.3pt}{\tiny \it #1}}}
\def\Tr{{\rm Tr}}

\def\str#1{\rule[#1mm]{0pt}{1mm}}
\def\Trr#1{{\rm Tr\str{-1.3}}_{\!R^{\mbox{\scriptsize
$(#1)$}}}}

\def\be{\begin{equation}}
\def\ee{\end{equation}}

\begin{document}

\makeatletter
\renewcommand{\theequation}{{\thesection}.{\arabic{equation}}}
\@addtoreset{equation}{section} \makeatother

\title{Braided Differential Operators on Quantum Algebras}

\author{
\rule{0pt}{7mm} Dimitri
Gurevich\thanks{gurevich@univ-valenciennes.fr}\\
{\small\it LAMAV, Universit\'e de Valenciennes,
59313 Valenciennes, France}\\
\rule{0pt}{7mm} Pavel Pyatov\thanks{pyatov@thsun1.jinr.ru}\\
{\small\it
Faculty of Mathematics, NRU HSE, 101000 Moscow, Russia~ \&}\\
{\small\it
Bogoliubov Laboratory of Theoretical Physics, JINR, 141980 Dubna, Russia}\\
\rule{0pt}{7mm} Pavel Saponov\thanks{Pavel.Saponov@ihep.ru}\\
{\small\it Division of Theoretical Physics, IHEP, 142284
Protvino, Russia} }
\date{}
\maketitle

\begin{abstract}
We propose a general scheme of constructing braided differential algebras via
algebras of "quantum exponentiated vector fields" and those of "quantum functions".
We treat a reflection equation algebra as a quantum analog of the algebra of
vector fields. The role of a quantum function algebra is
played by a general quantum matrix algebra. As an example we mention the
so-called RTT algebra of quantized functions on the linear matrix group
$GL(m)$. In this case our construction essentially coincides with the quantum differential algebra
introduced by S.~Woronowicz. If the role of a quantum function algebra is played by another
copy of the reflection equation algebra we get two different braided differential algebras. One of them
is defined via a quantum analog of  (co)adjoint vector fields, the other algebra is
defined via a quantum analog  of right-invariant vector fields. We show that the former algebra
can be identified with a subalgebra of the latter one. Also, we show that
"quantum adjoint vector fields" can be restricted to the so-called
"braided orbits" which are counterparts of generic  $GL(m)$-orbits in $gl^*(m)$.
Such braided orbits endowed with these
restricted vector fields constitute  a new class of braided differential algebras.
\end{abstract}

{\bf AMS Mathematics Subject Classification, 2010: 17B37, 81R50}

\section{Introduction}
\label{sec:i}

Since the creation of the quantum group theory plenty of different quantum
algebras related to R-matrices (i.e., solutions of the quantum Yang-Baxter equation)
have been introduced in the mathematical and physical literature. A remarkable family
of such algebras was introduced in \cite{STS} under the name of Heisenberg doubles.

As an associative algebra the Heisenberg double is generated by elements of two
dual Hopf algebras $H$ and $H^*$. In order to define an associative  product
on the space $H\ot H^*$ one needs a permutation operator
$$
H\ot H^* \to H^*\ot H
$$
transposing elements of two components. Such an operator can be defined via the pairing
$$
H\otimes H^*\rightarrow {\Bbb K}
$$
putting the algebras $H$ and $H^*$ in the duality ($\Bbb K$ is the ground field). Also,
assuming one of these algebras to be the quantized function algebra, namely, the famous
RTT algebra, one can extract from the dual object a space of "quantum exponentiated vector
fields" via a construction similar to that of \cite{FRT}.

A close approach to constructing a quantum version of differential calculus on a matrix pseudogroup
was  initiated by S.~Woronowicz \cite{W}. In that paper the central object --- quantum differential
algebra --- consists of three ingredients: a quantized function algebra, an algebra of quantum
exponentiated vector fields which is in fact the reflection equation (RE) algebra\footnote{This
algebra was studied by S.Majid under the name of  braided matrix algebra (see \cite{M} and the
references therein). The term "reflection equation algebra" was introduced by P.Kulish and his
coauthors (see \cite{KS}).} and an algebra of  "quantum differential forms".

We would also refer the reader to the papers \cite{FP,IP2} where a different approach to the quantum
calculus was suggested. In particular, it was shown that the classical Leibnitz rule for the external
differential must be modified.

In the present paper we disregard quantum differential forms and generalize the other
components of this calculus as follows. We always keep the RE algebra\footnote{The RE algebra
is parameterized by an R-matrix $R$. Below we are dealing with R-matrices of $GL(m)$-type
defined in Section \ref{sec:hr-qma}. However, a big part of our results can be generalized to algebras
associated with R-matrices of a more general form.} ${\cal L}(R)$ as an algebra of quantum
exponentiated vector fields. However, we introduce different candidates on the role of a quantum
function algebra ${\cal M}$, endowed with an appropriate action of the RE algebra. For any of such
a candidate the key point consists in constructing a permutation operator
\be
\mbox{\sf R}:\quad {\cal L}(R)\otimes {\cal M}\rightarrow {\cal M} \otimes  {\cal L}(R) \label{perm}
\ee
which enables us to endow the space ${\cal L}(R)\otimes {\cal M}$ with an associative product. We
denote the resulting associative algebra by ${\cal B}({\cal L}(R), {\cal M})$ and call it {\em a braided
differential (BD) algebra}.

We emphasize that the operator (\ref{perm}) also enables us to write down a sort of the Leibnitz rule for
elements of ${\cal L}(R)$. Having defined an action of these elements on generators of the algebra
${\cal M}$ and extending this action  to the higher components of ${\cal M}$ via such a "Leibnitz rule",
we get a representation of ${\cal L}(R)$ into the algebra ${\cal M}$. Note that in this construction we
have no need of a bialgebra structure (usual or braided) of the algebra ${\cal M}$. Instead, we use a
braided bialgebra structure of ${\cal L}(R)$ in order to apply its representation theory developed in
\cite{GPS1}\footnote{In \cite{GPS1} we constructed a representation category for the so-called modified
RE algebra which in fact coincides with the non-modified RE algebra but written in another basis. All
construction and results of \cite{GPS1} can be directly adapted to the non-modified form of the RE
algebra. Below we  refer to this paper without saying it each time.}.

Now, describe in more detail different types of the quantum algebras ${\cal M}$ we are dealing with.
First, we consider algebras ${\cal M}$ generated by the basic objects $V$ and $V^*$ of the ${\cal L}(R)$
representation category constructed in \cite{GPS1}. Note that the free tensor algebras $T(V)$ and $T(V^*)$,
as well as the "R-symmetric" and  "R-skew-symmetric" algebras of the space $V$ (or $V^*$), are examples
of such algebras ${\cal M}$.

Second, we consider quantum {\em matrix} algebras ${\cal M}$, each constructed via a pair of compatible
R-matrices (see definition (\ref{cc}) in Section \ref{sec:hr-qma}). As a particular case of such an algebra
${\cal M}$ we get the RTT algebra. In this case the resulting BD algebra coincides with the Heisenberg double
studied in \cite{IP}. If we take another copy of the RE algebra as ${\cal M}$ we get one more example of a BD
algebra. In this case two copies of the RE algebras are involved -- one of them (denoted ${\cal L}(R)$) plays the
role of quantum exponentiated vector fields, the other one (denoted ${\cal M}(R)$) plays the role of a quantum
function algebra.

The characteristic property of these two and other similar examples of BD algebras is that the elements of
${\cal L}(R)$ act on the quantum matrix algebra ${\cal M}$ on the left side and are in a sense analogs of
right-invariant exponentiated vector fields. In what follows such BD algebras are denoted ${\cal B}_r({\cal L}(R),
{\cal M})$ where the subscript $r$ means "right-invariant".

However, if the quantum matrix algebra ${\cal M}$ is just the RE algebra ${\cal M}(R)$ we can define another
action of ${\cal L}(R)$ on ${\cal M}(R)$, namely the ``adjoint'' action, which is  an analog of the usual adjoint action
of one copy of $gl(m)$ onto another one. The corresponding BD algebra  is denoted ${\cal B}_{ad}({\cal L}(R),
{\cal M}(R))$. Thus, we have two versions of the BD algebra composed from the algebras ${\cal L}(R)$ and
${\cal M}(R)$. One of them ${\cal B}_r({\cal L}(R), {\cal M}(R))$ is based on the right-invariant action, the other
one --- ${\cal B}_{ad}({\cal L}(R), {\cal M}(R))$ --- on the adjoint action. We show that the algebra
${\cal B}_{ad}({\cal L}(R), {\cal M}(R))$ can be embedded into the properly extended algebra
${\cal B}_r({\cal L}(R), {\cal M}(R))$ as a subalgebra.

Similarly to the classical case the adjoint action of the algebra ${\cal L}(R)$ onto ${\cal M}(R)$ preserves central
elements of the latter algebra. Using this fact it is possible to reduce the quantum adjoint vector fields to "braided
orbits", i.e. the quotients of  ${\cal M}(R)$ which are quantum counterparts of generic $GL(m)$ orbits in $gl^*(m)$
(see  \cite{GS1,GS3} for detail). Thus, we get one more family of BD algebras, in which the role of quantum function
algebras is played by "braided orbits".

The paper is organized as follows. In the next Section we recall some elements of "braided geometry" as presented
in \cite{GS2}. In particular, we exhibit a regular way of constructing a quantum matrix algebra via a pair of compatible
R-matrices. In Section \ref{sec:rea-rep} we concentrate ourselves on properties of the RE algebra including its
representation category. This enables us to construct some examples of the BD algebras considered in Section
\ref{sec:qd-rea}. In Section \ref{sec:5} we present a construction of the BD algebra over a general quantum matrix
algebra and exhibit the mentioned relation between BD algebra based on two copies of the RE algebra but equipped
with different types of action. We complete the paper with an example of BD algebra on a quantum hyperboloid.

\medskip

\noindent
{\bf Acknowledgement.}
This work was jointly supported by CNRS and RFBR grants GDRI-471, 09-01-93107-NCNIL-a and 08-01-00392-a.
The work of PP was supported by RFBR (grant No.11-01-00980-a) and by NRU HSE (personal research grant
No.~10-01-0013). The work of PS was supported by the joint RFBR and DFG grant 08-01-91953-a. Also PP and PS
are grateful to the Valenciennes University (France) where this paper was partially written for the warm hospitality and
support.

\section{$GL(m)$-type $R$-matrices and quantum matrix algebras}
\label{sec:hr-qma}

In this section we give a short list of definitions and notation to be used
below. More details and proofs can be found in the cited literature.

Let $\Bbb K$ denote the field of complex or real numbers and $V$ be a finite
dimensional vector space over the field $\Bbb K$:  $\dim_{\K}V=N$. Given a
linear operator $X\in \End(V^{\otimes \,k})$, $\forall k\geqslant 1$, we extend it
up to different operators belonging to $\End(V^{\otimes (k+p)})$, $p\geqslant 0$,
in a natural way
\be
X_{i\dots i+k-1}=I_V^{\otimes (i-1)}\otimes X\otimes I_V^{\otimes (p-i+1)},
\quad 1\leqslant i \leqslant p+1,
\label{matr-conv}
\ee
where $I_V$ stands for the identity operator on $V$. In what follows
we shall abbreviate $I_V$ to $I$ and simplify $X_{i\,i+1}$ to $X_i$ for
$X\in\End(\vv)$. Hereafter, all tensor products are taken over the ground field $\Bbb{K}$.

An invertible operator $R\in \mathrm{Aut}(V^{\otimes 2})$ is called an
{\it R-matrix} if it satisfies the {\it Yang-Baxter equation} in $\End(V^{\otimes 3})$
\be
R_1R_2R_1 - R_2R_1R_2 = 0.
\label{YBE}
\ee

In the present paper we are dealing with {\it Hecke type} R-matrices which
obey the quadratic Hecke condition
\be
(R - q\, I^{\otimes 2})(R + q^{-1}\,I^{\otimes 2})=0\,, \quad q\in {\Bbb K}
\setminus 0\,.
\label{Hec}
\ee
A numerical parameter $q$ is assumed to be {\it generic}, that is either $q=1$
or $q^k\not=1$, $\forall \,k\in {\Bbb N}$. In particular, for a generic value of the
parameter the $q$-analogs of integers
$$
k_q=q^{k-1}+q^{k-3}+\dots+ q^{1-k} = \frac{q^k-q^{-k}}{q-q^{-1}}
$$
are non-zero for any integer $k\in \Bbb Z$.

An example of the Hecke type R-matrix for $q=1$ is given by the flip
(transposition operator):
\be
P: V^{\otimes 2}\rightarrow V^{\otimes 2},\quad
P(v_1\otimes v_2) = v_2\otimes v_1.
\label{flip-P}
\ee
A well-known example for $q\not = 1$ is the Drinfeld-Jimbo R-matrix
\be
R(q) = q\,\sum_{i=1}^NE_{ii}\otimes E_{ii} +\sum_{i\not=j}^NE_{ij}\otimes E_{ji}
+(q-q^{-1})\sum_{1\leqslant i<j\leqslant N}E_{ii}\otimes E_{jj},
\label{DJ-R}
\ee
where $E_{ij}\in {\rm Mat}_N({\Bbb K})$ are the standard matrix units.
Note that $R(q)$ is a continuous matrix function in $q$
and $\lim_{q\rightarrow 1}R(q) = P$. There are known other Hecke type
R-matrices which are continuous matrix functions in $q$ and turn into the
flip $P$ at the limit $q\rightarrow 1$. All such Hecke type R-matrices will be
referred to as deformations of the flip $P$.

The Hecke type R-matrices are closely connected with the representation theory
of the $A_{n-1}$ series Hecke algebras ${\cal H}_n(q)$, $n\geqslant 2$. Recall,
that the Hecke algebra ${\cal H}_n(q)$ is the quotient of the group algebra
${\Bbb K}[{\cal B}_n]$ of the braid group ${\cal B}_n$, $n\geqslant 2$,
$$
{\cal B}_n = \langle \{\sigma^{\pm}_i\}_{1\leqslant i\leqslant n-1}: \sigma_i
\sigma_{i+1}\sigma_i =  \sigma_{i+1}\sigma_i\sigma_{i+1},\quad\sigma_i
\sigma_j = \sigma_j\sigma_i\quad i\not=j\pm 1,\quad \sigma_i^{\pm 1}
\sigma_i^{\mp 1} = 1_{\mbox{\tiny ${\cal B}$}}\rangle
$$
over the two sided ideal generated by the elements
$$
\sigma_i^{-1} - \sigma_i + (q-q^{-1})1_{\mbox{\tiny ${\cal B}$}}
$$
where $1_{\mbox{\tiny ${\cal B}$}}$ stands for the unit element of the braid group.

At a generic value of $q$ the algebra ${\cal H}_n(q)$ is known to be semisimple
and isomorphic to the group algebra ${\Bbb K}[{\frak S}_n]$ of the $n$-th order
permutation group. As a consequence, the primitive idempotents $e_\lambda^a:
e_\lambda^ae_\mu^b = \delta_{\lambda\mu}\delta^{ab}e_{\lambda}^a$
of the Hecke algebra ${\cal H}_n(q)$ are labeled by partitions $\lambda\vdash n$
and by an integer index $a: 1\leqslant a\leqslant d_\lambda$ where
$d_\lambda$ equals to the number of standard Young tableaux corresponding
to the partition $\lambda$.

Any Hecke type R-matrix $R$ realizes a so-called {\it local R-matrix representation}
${\rho\str{-1.3}}_{\! R}$ of a Hecke algebra ${\cal H}_n(q)$ by the following rule:
\be
{\rho\str{-1.3}}_{\! R}: {\cal H}_n(q)\rightarrow \End(V^{\otimes n}),\qquad
{\rho\str{-1.3}}_{\! R}(\sigma_i) = R_i,\quad 1\leqslant i\leqslant n-1.
\label{r-rep}
\ee
The detailed treatment of the Hecke algebra and its representations with a
list of original papers can be found, e.g., in the review \cite{OP}.

We constrain ourselves to considering a subfamily of the Hecke type R-matrices
--- so-called   R-matrices of the $GL(m)$-type.

A Hecke type R-matrix $R$ is said to be of the $GL(m)$-type if the operator
$$
A_m(R) = {\rho\str{-1.3}}_{\! R}(e_{(1^m)}(\sigma))
$$
is a rank one projector in ${\rm End}(V^{\otimes m})$, while
$A_{m+1}(R)\in {\rm End}(V^{\otimes (m+1)})$ is the zero operator (explicit
formulae for $A_k(R)$, $k\geqslant 1$, can be found, for example, in \cite{OP}). Note,
that the $N^2\times N^2$ Drinfeld-Jimbo R-matrix (\ref{DJ-R}) is of the $GL(N)$-type,
but for general $GL(m)$-type R-matrix we have $m\leqslant N$ (see \cite{G} for examples).

Also, we need the following definition.
An operator $R\in \End(V^{\otimes 2})$ is called {\it skew-invertible} provided
that there exists an operator $\Psi_{\sscr{R}}\in \End(V^{\otimes 2})$ such that
\be
{\rm Tr}_{(2)}R_{12}{\Psi_{\sscr{R}}}_{23} = P_{13} =
{\rm Tr}_{(2)}{\Psi_{\sscr{R}}}_{12} R_{23}\, ,
\label{def:psi}
\ee
where the subscript in the notation of the trace indicates the factor in the tensor
product $V^{\otimes 3}$, where the trace operation is applied.

The operators
\be
B_\sscr{R} = \Tr_{(1)}{\Psi_\sscr{R}}_{12},\qquad C_\sscr{R}  = \Tr_{(2)}
{\Psi_\sscr{R}}_{12},
\label{def:BC}
\ee
possess the property
\be
\Tr_{(1)}{B_\sscr{R}}_1R_{12} = I =\Tr_{(2)}{C_\sscr{R}}_2R_{12}.
\label{BRC}
\ee
A skew-invertible R-matrix is called {\it strictly} skew-invertible if the corresponding
operator $B_\sscr{R}$ (or $C_\sscr{R}$) is invertible.

As was shown in  \cite{G}, any $GL(m)$-type R-matrix is skew-invertible. Moreover,
the operators $B_\sscr{R}$ and $C_\sscr{R}$ are subject to the relation \cite{GPS1}
$$
B_\sscr{R}\cdot C_\sscr{R} = q^{-2m}I
$$
Consequently, any $GL(m)$-type R-matrix is strictly skew-invertible.

The operators $B_\sscr{R}$ and $C_\sscr{R}$ play an important role in the theory
of quantum matrix algebras considered below. In particular, the operator $C_\sscr{R}$
appears in definition of the {\it R-trace} ${\Tr\str{-1.3}}_{\!R}$:
\be
{\Tr\str{-1.3}}_{\!R}: {\rm Mat}_N(A) \rightarrow A, \qquad
{\Tr\str{-1.3}}_{\!R}(X) = \Tr(C_\sscr{R}X),\quad
\forall\,X\in {\rm Mat}_N(A),
\ee
where $A$ is a vector space over the field $\K$ and
${\rm Mat}_N(A) = {\rm Mat}_N(\K)\otimes_\K A$.

Turn now to the definition of a quantum matrix algebra \cite{IOP}. An {\it ordered} pair
$\{R,F\}$ of two R-matrices $R,F\in {\rm Aut}(V^{\otimes 2})$ is called {\it
compatible} if they satisfy the following relations (compatibility conditions):
\be
R_1F_2F_1 = F_2F_1R_2\, ,
\qquad
F_1F_2R_1 = R_2F_1F_2\, .
\label{cc}
\ee

\begin{definition}
\label{def:QMA}
{\bf \cite{IOP}} \hspace{3pt}{\rm
Given a compatible pair $\{R,F\}$ of strictly skew-invertible R-matrices
$R,F\in {\rm Aut}(V^{\otimes 2})$, $\dim_{\Bbb K}V = N$, the {\em quantum matrix (QM)
algebra} ${\cal M}(R,F)$ is a unital associative algebra generated by a unit
element $1_\sscr{${\cal M}$}$ and by $N^2$ entries of the matrix
$M=\|M_i^{\,j}\|_{1\leqslant  i,j\leqslant N}$ subject to the relations
\be
\label{def:QM}
R_1M_{\overline 1}M_{\overline 2} -M_{\overline 1}
M_{\overline 2} R_1 = 0\, ,
\ee
where we use the notation
\be
M_{\overline 1} = M_1,  \qquad M_{\overline {k+1}} =
F_{k} M_{\overline {k}}\, F^{-1}_{k},\quad k\geqslant 1,
\label{copy}
\ee
for the ``copies'' of the matrix $M$. In what follows the matrix $M$ will be
called a generating matrix of the QM algebra ${\cal M}(R,F)$.
(In the above formulae $R_{1}$ and $F_{k}$, $k\geqslant 1$, are treated in the sense of formula
(\ref{matr-conv}).)

 Also, observe that by fixing a basis $\{x_i\}_{1\leqslant i\leqslant N}$  in the space $V$ and the corresponding
basis  $\{x_i\otimes x_j\}$ in that $V^{\ot 2}$ we can represent the operators $R,F\in {\rm Aut}(\vv)$ by
{\it numerical} matrices $\|R_{ij}^{\,rs}\|$ and
$\|F_{ij}^{\,rs}\|$ where
\be
R(x_i\otimes x_j) = R_{ij}^{\,rs}\,x_r\otimes x_s,\qquad
F(x_i\otimes x_j) = F_{ij}^{\,rs}\,x_r\otimes x_s.
\label{r-bas}
\ee
}
\end{definition}

The defining relations (\ref{def:QM}) and compatibility conditions (\ref{cc}) imply the same type
relations for consecutive pairs of the copies of $M$ \cite{IOP}
\be
R_k\,M_{\overline k}M_{\overline{k+1}} -
M_{\overline k}M_{\overline{k+1}}\, R_k = 0\,.
\nonumber
\ee

It follows from Definition \ref{def:QMA} that any QM algebra ${\cal M}(R,F)$ is a finitely
generated quadratic (in the generators $M_i^{\,j}$) graded algebra  and consequently, it can be presented
as a sum of homogeneous components
$$
{\cal M}(R,F) = \sum_{p\geqslant 0} {\cal M}^p(R,F),\qquad
{\cal M}^0(R,F)\cong {\Bbb K}.
$$

The widely known example of QM algebra is the quantized algebra of functions on the
matrix algebra $M_N(\Bbb K)$ \cite{FRT}. This algebra is associated with a compatible
pair $\{R,P\}$, $P$ being a flip (\ref{flip-P}). Let $T = \|T_i^{\,j}\|_1^N$ be the generating
matrix of ${\cal M}(R,P)$. Then the relations (\ref{def:QM}) take the form
\be
R_1T_1T_2 - T_1T_2R_1=0,
\label{RTT}
\ee
since  $T_{\overline 2}  = P_{12}T_1P_{12} = T_2$. We call the QM algebra ${\cal M}(R,P)$ the
RTT algebra and denote it by ${\cal T}(R)$. At any choice of an R-matrix $R$ the RTT algebra
${\cal T}(R)$ is a bialgebra with the coproduct $\Delta$ and counit $\varepsilon$
\be
\Delta(T) = T\stackrel{.}{\otimes}T,\qquad \varepsilon(T) = I.
\label{bi-rtt}
\ee
The symbol $\stackrel{.}{\otimes}$ stands for the following operation
\be
(A\stackrel{.}{\otimes}B)_{i}^{\,j} = \sum_{k=1}^NA_{i}^{\,k}\otimes B_k^{\,j}
\label{tens-dot}
\ee
where $A$ and $B$ are arbitrary $N_1\times N$ and $N\times N_2$ matrices
respectively.

Note that for any R-matrix of $GL(m)$-type a quantum determinant ${\rm det}_q\,T$ can be defined
(see \cite{G}). This enables us to introduce an antipode and a Hopf algebra structure in the algebra
${\cal T}(R)$ extended by $({\rm det}_q\,T)^{-1}$. The resulting algebra is the most popular quantum
analog of the function algebra on  $GL(m)$ (see  \cite{FRT}).

Below we consider a family of the QM algebras associated with a compatible pair
$\{R,F\}$ where $R$ is a $GL(m)$-type R-matrix. The corresponding algebra
${\cal M}(R,F)$ is referred to as a $GL(m)$-type QM algebra.

A useful tool for studying the structure of the $GL(m)$-type QM algebra is the {\it
characteristic subalgebra} ${\rm Char }({\cal M})\subset {\cal M}(R,F)$ \cite{IOP}.
By definition, this is a linear span of the unit element and the following elements
$$
x(h_k) = \Trr{1\dots k}(M_{\overline 1}\dots M_{\overline
k}\,{\rho\str{-1.3}}_{\! R}(h_k)),\quad k \in {\Bbb N},
$$
where $h_k$ runs over all elements of the Hecke algebra ${\cal H}_k(q)$ and
${\rho\str{-1.3}}_{\! R}$ is the R-matrix representation (\ref{r-rep}) of ${\cal H}_k(q)
$
in $\End(V^{\otimes k})$.

Among all elements of the characteristic subalgebra we distinguish the following
families:
\begin{itemize}
\item
The {\it power sums} of the quantum matrix
\be
p_0(M) = 1_\sscr{${\cal M}$}{\Tr\str{-1.3}}_{\!R}(I),\quad p_k(M)=\Trr{1\dots
k}
(M_1M_{\bar 2}\dots M_{\bar k} {\rho\str{-1.3}}_{\! R}(\sigma_{k-1}\dots \sigma_2
\sigma_1)), \quad k\geqslant 1;
\label{pow-sums}
\ee
\item
The {\it elementary symmetric functions}
\be
a_0(M) = 1_\sscr{$\cal M$},
\quad a_k(M) = \Trr{1\dots k}(M_1M_{\bar 2}\dots M_{\bar k}{\rho\str{-1.3}}
_{\! R}(e_{(1^k)})),
\quad 1\leqslant k\leqslant m.
\label{elem-f}
\ee

\end{itemize}

The main properties of a $GL(m)$-type QM algebra to be used below are collected in
the following proposition.
\begin{proposition} {\bf \cite{IOP}}
Let ${\cal M}(R,F)$ be a $GL(m)$-type quantum matrix algebra. Then the
following statements hold true.
\begin{enumerate}
\item
The characteristic subalgebra ${\rm Char}({\cal M})$ is abelian.
\item
The characteristic subalgebra is generated by the set of {\it power sums}
$\{p_k(M)\}_{0\leqslant k\leqslant m}$ or, equivalently, by the set of elementary
symmetric functions $\{a_k(M)\}_{0\leqslant k\leqslant m}$.
\item
The power sums are related with elementary symmetric functions by the quantum
Newton identities
\be
(-1)^{k-1}k_qa_k(M) = \sum_{i=0}^{k-1}(-q)^ip_{k-i}(M)a_i(M),\quad 1\leqslant
k\leqslant m.
\label{q-newt}
\ee
\item
The generating matrix $M$ satisfies the Cayley-Hamilton identity
\be
\sum_{k=0}^m(-q)^kM^{\overline{m-k}}a_k(M) = 0,
\label{qma-ch}
\ee
where
$$
M^{\overline 0} = 1_\sscr{${\cal M}$}I,\quad
M^{\overline k} = \Trr{2\dots k}(M_{\overline 1}M_{\overline 2}\dots M_{\overline k}
{\rho\str{-1.3}}_{\! R}(\sigma_{k-1}\sigma_{k-2}\dots \sigma_1)).
$$
\end{enumerate}
\end{proposition}

\section{RE algebra and its representation theory}
\label{sec:rea-rep}

In this section we give a short review of a particular case of the QM algebra --- the
RE algebra \cite{KS,M}. We discuss its structure and the representation theory. Our
exposition will mainly follow the paper \cite{GPS1} where the reader can find detailed
proofs and further references to the literature on the RE algebra.

The $GL(m)$-type RE algebra ${\cal L}(R)$ is associated with the compatible pair
$\{R,R\}$, where $R\in {\rm End}(V^{\otimes 2})$ is a $GL(m)$-type R-matrix (recall,
that $\dim_{\Bbb K}V = N\geqslant m$). Denoting the matrix of the generators of the RE algebra by
$L=\|L_i^{\,j}\|_1^N$ we rewrite the general commutation relations
(\ref{def:QM}) in the equivalent form:
\be
R_1L_1R_1L_1 - L_1R_1L_1R_1 = 0.
\label{REA}
\ee
Note, that the algebra ${\cal L}(R)$ has the structure of the left coadjoint comodule
over the RTT algebra ${\cal T}(R)$ defined by (\ref{RTT}). On the first order
homogeneous component ${\cal L}^1(R)$ (the linear span of the generators of the RE algebra)
the coaction  $\delta_\ell: {\cal L}(R) \rightarrow {\cal T}(R)\otimes {\cal L}(R)$
reads
\be
\delta_\ell(L_i^{\,j}) = \sum_{k,p=1}^NT_i^{\,k}S(T_p^{\,j})\otimes L_k^{\,p},
\label{rtt-rea-coac}
\ee
where $S(T)$ stands for the antipodal mapping of the RTT Hopf algebra.

The main peculiarities of the RE algebra (comparing with a general QM algebra) are listed below.
\begin{enumerate}
\item
The quantum matrix powers $L^{\overline k}$, $k\geqslant 1$,  (\ref{qma-ch}) are
simplified to the usual matrix products $L^{\overline k}=L^k = L\cdot L^{k-1}$, where
$L^0 = 1_\sscr{${\cal L}$}I$. Consequently, the power sums $p_k(L)$ (\ref{pow-sums})
take the form $p_k(L) =  {\Tr\str{-1.3}}_{\!R}(L^k)$, $k\geqslant 0$.
\item
The power sums $p_k(L)$, $k\geqslant 0$ are {\it central} elements in ${\cal L}
(R)$ \cite{KS}. As a consequence, the abelian characteristic subalgebra
${\rm Char}({\cal L})$ is {\it central} in the RE algebra ${\cal L}(R)$.
\item
Due to the property 1 the Cayley-Hamilton identity (\ref{qma-ch}) for the matrix $L$ takes the usual form:
\be
\sum_{i=0}^m (-q)^ia_i(L)L^{m-i} =0.
\label{CH-id}
\ee
\end{enumerate}

For the $GL(m)$-type RE algebra we introduce $m$ spectral values $\mu_i$,
$1\leqslant i \leqslant m$, of the quantum matrix $L$ considered as elements
of a central extension of the ${\rm Char}({\cal L})$ (see \cite{GPS2} for more
detail and a generalization to the case of $GL(m|n)$-type R-matrix). These spectral
values are defined by the following system of polynomial relations
\be
a_k(L) = q^{-k}\sum_{1\leqslant i_1<i_2<\dots i_k\leqslant m}\mu_{i_1}\mu_{i_2}
\dots \mu_{i_k},\quad 1\leqslant k\leqslant m.
\label{elem-spectr}
\ee

Then, any element of the characteristic subalgebra can be parameterized by a
symmetric polynomial in spectral values. In particular, the parameterization of
power sums $p_k(L)$ reads \cite{GS1}:
$$
p_k(L) = \sum_{i=1}^m d_i\mu_i^k,\qquad d_i = q^{-m}\prod_{j\not=i}
^m\frac{q\mu_i-q^{-1}\mu_j} {\mu_i-\mu_j}.
$$
Besides, (\ref{CH-id}) and (\ref{elem-spectr}) allow the Cayley-Hamilton identity (\ref{CH-id})
to be written in a factorized form
$$
\prod_{i=1}^m(L-\mu_iI) = 0.
$$

The next property of the RE algebra ${\cal L}(R)$ allows us to construct a category of its finite
dimensional representations. Namely, the RE algebra has a structure of a {\it braided bialgebra}
\cite{M,GPS1}. To define this structure we need some more notation.

Introduce a finite dimensional vector space $W(L)$:
\be
W(L) = {\rm span}_{\Bbb K}\{L_i^{\,j}, 1\leqslant i,j\leqslant N\},
\qquad \dim_{\Bbb K}W(L) = N^2,
\label{def:wk}
\ee
and consider the free tensor algebra $TW(L)$ generated by the space $W(L)$.
In each homogeneous component $W(L)^{\otimes k}\subset TW(L)$,
$k\geqslant 1$, we take the basis, formed by entries of the $N^k\times N^k$
matrix $L_{1\rightarrow\overline k}$:
\be
L_{1\rightarrow\overline k}\stackrel{\mbox{\tiny def}}{=}
L_1\stackrel{.}{\otimes}L_{\overline 2}\stackrel{.}{\otimes}\dots\stackrel{.}
{\otimes} L_{\overline {k-1}}\stackrel{.}{\otimes}
L_{\overline k},\qquad L_{\overline {i+1}} = R_iL_{\overline i}R_i^{-1}.
\label{bas-ov}
\ee
Another possible choice of the basis set in $W(L)^{\otimes k}$ is given by
elements of the matrix
\be
L_{\underline k\rightarrow 1}\stackrel{\mbox{\tiny def}}{=}
L_{\underline k}\stackrel{.}{\otimes}L_{\underline{k-1}}\stackrel{.}{\otimes}
\dots\stackrel{.}{\otimes} L_{\underline 2}\stackrel{.}{\otimes}L_1,\qquad
L_{\underline {i+1}} = R^{-1}_iL_{\underline i}R_i.
\label{bas-un}
\ee

The RE algebra ${\cal L}(R)$ is isomorphic to the quotient
$$
{\cal L}(R)\cong TW(L)/\langle R_1L_{1\rightarrow \overline 2} -
L_{1\rightarrow \overline 2}R_1\rangle,
$$
where $\langle J\rangle $ denotes the two-sided ideal in the tensor algebra
$TW(L)$ generated by a subset $J\subset TW(L)$. Being projected from
$TW(L)$ to the algebra ${\cal L}(R)$, the sets (\ref{bas-ov}) and (\ref{bas-un})
have coinciding images
$$
L_1L_{\overline 2}\dots L_{\overline {k-1}}L_{\overline k} = L_{\underline k}
L_{\underline{k-1}}\dots L_{\underline 2} L_1.
$$
This can be proved with the use of permutation relations (\ref{REA}) and
Yang-Baxter equation for $R$.  The homogeneous component ${\cal L}^k(R)$
is a linear span of matrix elements of the matrix
$$
L_{1\rightarrow {\overline k}} = L_1L_{\overline 2}\dots L_{\overline k}
$$
where we keep the same notation as in (\ref{bas-ov}).

The braided bialgebra structure of the RE algebra ${\cal L}(R)$ is defined by two
homomorphic maps: the {\it coproduct} $\Delta: {\cal L}(R)\rightarrow
\mbox{\bf L}(R)$ and the {\it counit} $\varepsilon: {\cal L}(R)\rightarrow \Bbb K$.
Here an associative unital algebra $\mbox{\bf L}(R)$ has the following structure
\cite{GPS1}:
\begin{enumerate}
\item
As a vector space over the field $\Bbb K$ the algebra $\mbox{\bf L}(R)$ is
isomorphic to the tensor product of two copies of the RE algebra
$$
\mbox{\bf L}(R)\cong {\cal L}(R)\otimes {\cal L}(R);
$$
\item
The algebra $\mbox{\bf L}(R)$ is endowed with a vector space automorphism
$$
\mbox{\sf R}: \;{\cal L}^1(R)\otimes {\cal L}^1(R)\rightarrow
{\cal L}^1(R)\otimes {\cal L}^1(R)
$$
where ${\cal L}^k(R)$ stands for the homogeneous component of degree $k$.
On the basis elements of the space ${\cal L}^1(R)\otimes {\cal L}^1(R)$ the automorphism
is defined by the rule
\be
\mbox{\sf R}(L_1\stackrel{.}{\otimes} L_{\overline 2}) =
L_{\overline 2}\stackrel{.}{\otimes} L_1.
\ee
It can be extended to the set of automorphisms
$$
\mbox{\sf R}_{\mbox{\sf\footnotesize(k,p)}}:\;
{\cal L}^k(R)\otimes {\cal L}^p(R)\rightarrow
{\cal L}^p(R)\otimes {\cal L}^k(R), \quad {k\geqslant0,p\geqslant 0}.
$$
by the following relations
\be
\mbox{\sf R}_{\mbox{\sf\footnotesize(k,p)}}(L_{1\rightarrow {\overline k}}
\stackrel{.}{\otimes} L_{\overline{(k+1)}\rightarrow\overline{(k+p)}}) =
L_{\overline{(k+1)}\rightarrow\overline {(k+p)}}\stackrel{.}{\otimes}
L_{1\rightarrow {\overline k}}.
\label{R-SW}
\ee
The subspace ${\cal L}^0(R)$ generated by the unit element $1_\sscr{$\cal L$}$ commutes
with any ${\cal L}^k(R)$.
\item
Let $a,d\in {\cal L}(R)$ be arbitrary elements of the RE algebra and  $b\in {\cal L}^k(R)$,
$c\in{\cal L}^p(R)$ be arbitrary elements of homogeneous components of ${\cal L}(R)$.
Then by definition the product $(a\otimes b)*(c\otimes d)$ is given by the rule
\be
(a\otimes b)*(c\otimes d) = ac_{(1)}\otimes  b_{(2)}d,
\label{br-prod}
\ee
where $a c_{(1)}$ and $b_{(2)}d$ are products of elements of the RE algebra, while
$c_{(1)}$ and $b_{(2)}$ are the Sweedler's notation for the image
of the automorphism $\mbox{\sf R}_{\mbox{\sf\footnotesize(k,p)}}$:
$$
\mbox{\sf R}_{\mbox{\sf\footnotesize(k,p)}}(b\otimes c) =  c_{(1)}\otimes b_{(2)}=
\sum_{i} c_i\otimes b_i,\quad c_i\in{\cal L}^p(R),\,b_i\in {\cal L}^k(R).
$$
\end{enumerate}
The product of arbitrary elements of the algebra $\mbox{\bf L}(R)$ is obtained
from the above definition by linearity.

The following proposition holds true \cite{GPS1} (in an equivalent form the proposition
was proved in \cite{M}).
\begin{proposition}
\label{pro:5}
Consider two linear maps $\Delta: {\cal L}(R)\rightarrow \mbox{\bf L}(R)$
and $\varepsilon: {\cal L}(R)\rightarrow \Bbb K$ defined by the relations
\be
\Delta(1_\sscr{$\cal L$}) = 1_\sscr{$\cal L$}\otimes
1_\sscr{$\cal L$},\qquad
\Delta(L_{1\rightarrow \overline{k}}) = L_{1\rightarrow \overline{k}}\stackrel{.}
{\otimes} L_{1\rightarrow \overline{k}},\quad k\geqslant 1
\label{copr}
\ee
and
\be
\varepsilon(1_\sscr{$\cal L$}) = 1,\qquad \varepsilon
(L_{1\rightarrow \overline{k}}) = I_{12\dots
k}.
\label{coun}
\ee
Then the maps $\Delta$ and $\varepsilon$ are homomorphisms of associative
unital algebras and they define a braided bialgebra structure on the RE algebra ${\cal L}
(R)$ with the coproduct $\Delta $ and the counit $\varepsilon$.
\end{proposition}

The representation theory of the $GL(m)$-type RE algebra ${\cal L}(R)$ generated by
$N^2$ generators $L_i^{\,j}$ can be developed in a monoidal rigid quasitensor (provided $q\not=1$)
category (see \cite{CP} for terminology) generated by an $N$-dimensional vector
space $V$, $R\in {\rm Aut}(\vv)$ (\cite{GLS}, see also \cite{GPS1} for generalization
to $GL(m|n)$ case). Following \cite{GPS1} we shall call this category the Schur-Weyl
category.

The term ``quasitensor'' means, that for any couple of objects $U_1$ and $U_2$ of
the category the functorial commutativity (iso)morphism $\mbox{\sf R}_{(U_1,U_2)}:
U_1\otimes U_2\rightarrow U_2\otimes U_1$ is {\it not} involutive (unless one of the
objects is the field $\Bbb K$). The above mappings
$\mbox{\sf R}_{\mbox{\sf\footnotesize(k,p)}}$ give examples of such morphisms.

The rigidity means that for any object $U$ its dual $U^*$ is also an object of the category
and moreover, there exist a left $\langle U^*\otimes U\rangle_l\rightarrow {\Bbb K}$ and a
right $\langle U\otimes U^*\rangle_r\rightarrow {\Bbb K}$ pairings which are morphisms of
the category (evaluation morphisms). Besides, there exist embeddings of the field
${\Bbb K}\rightarrow U\otimes U^*$ and ${\Bbb K}\rightarrow U^*\otimes U$ which are also
morphisms (co-evaluation morphisms).

Given a basis $\{x_i\}_{1\leqslant i\leqslant N}$ (\ref{r-bas}) in the generating space $V$,
then a basis $\{y^i\}_{1\leqslant i\leqslant N}$ in the dual space $V^*$ can be chosen in
such a way that
\be
\langle x_i,y^j\rangle_r =\delta_i^{\,j},\qquad \langle y^i,x_j\rangle_l =(B_\sscr{R}) \label{duba}
^{\,i}_j.
\ee
The aforementioned co-evaluation morphisms in these basis sets read
$$
{\Bbb K}\rightarrow V\otimes V^*:\; 1\rightarrow (C_\sscr{R})_i^{\,j}x_j\otimes y^i,
\qquad {\Bbb K}\rightarrow V^*\otimes V: \; 1\rightarrow y^i\otimes x_i.
$$

In \cite{GPS1} it was argued that the space $W(L)$ defined in (\ref{def:wk}) can
be treated as an object of the Schur-Weyl category isomorphic to $V\otimes V^*$.
This fact allows us to construct the categorical commutativity morphisms
$\mbox{\sf R}_{\,(W(L),V^{\otimes p})}$ which play the crucial role
in extending the ${\cal L}(R)$-module structure from the space $V$ to any its tensor
power. In particular, as minimal ``building blocks'' we have the following relations
\begin{eqnarray}
&&\mbox{\sf R}_{\,(W(L),V)}(L_{\underline 2}\stackrel{.}{\otimes} x_1) = x_1\otimes L_2,
\label{K-V}\\
&&\mbox{\sf R}_{\,(W(L),V^*)}(L_2\otimes y_1) = y_1\stackrel{.}{\otimes}
L_{\underline 2}.\label{L-Vd}
\end{eqnarray}
For two copies of the space $W(L)$ we get \cite{GPS1} (compare with (\ref{R-SW}))
\be
\mbox{\sf R}_{\,(W(L),W(L))}(L_1\stackrel{.}{\otimes} L_{\overline 2}) =
L_{\overline 2}\stackrel{.}{\otimes} L_1.
\label{comm-KM}
\ee

Also, the isomorphism $W(L)\cong V\otimes V^*$ enables us to define the
{\it adjoint} representation of the RE algebra ${\cal L}(R)$ on the space $W(L)$. Besides, the adjoint
action sends to zero the ideal ${\cal J}$ generated by the left-hand side of relations (\ref{REA}).
So, the adjoint action, being extended to the whole algebra ${\cal L}(R)\cong TW(L)/{\cal J}$
via the braided coproduct (\ref{copr}) and the commutativity morphism (\ref{comm-KM}),
respects the algebraic structure of the RE algebra.

Now, we present some explicit formulae of the ${\cal L}(R)$ representations in various spaces.
In the basis $\{x_i\}_{1\leqslant i\leqslant N}$ (\ref{r-bas})  the
action $\triangleright $ of the linear operator  $L_i^{\,j}$ is defined as follows
$$
L_i^{\,j}\triangleright x_p = \delta_i^{\,j}x_p - (q-q^{-1})(B_\sscr{R})_p^{\,j} x_i
$$
where $(B_\sscr{R})_p^{\,j}=\sum_a(\Psi_\sscr{R})_{ap}^{\,aj}$ according to
(\ref{def:BC}). Using the property (\ref{BRC}) one can easily show that the above
action provides the space $V$ with the left ${\cal L}(R)$-module structure. We
rewrite the above action in an equivalent covariant matrix form
\be
L_1R_1\triangleright x_1 = R_1^{-1}x_1.
\label{b-rep}
\ee
The compact formula (\ref{b-rep}) is a concise notation for the following
expression
$$
\sum_{a,b=1}^N(L_{i_1}^{\,a}R_{ai_2}^{\,bj_2})\triangleright x_b =
\sum_{a=1}^N(R^{-1})_{i_1i_2}^{\,aj_2}x_a.
$$

The representation (\ref{b-rep}) is irreducible provided that the matrix $B_\sscr{R}$
is nonsingular.

\begin{remark}{\rm
The RE algebra ${\cal L}(R)$ is defined by the quadratic relations (\ref{REA}), so it
admits an evident rescaling automorphism $L\mapsto \eta L$, with arbitrary
non-zero $\eta\in {\Bbb K}$. As a consequence, the action
\be
L_1R_1\triangleright x_1 = \eta\,R_1^{-1} x_1, \quad \eta\in {\Bbb K}\setminus 0
\label{rep-gen}
\ee
is also a representation of the algebra ${\cal L}(R)$.}
\end{remark}

To extend the ${\cal L}(R)$-module structure to $V^{\otimes p}$, $p\geqslant 2$,
we use the coproduct operation (\ref{copr}) and an inductive procedure. Let spaces
$U$ and $W$ be left ${\cal L}(R)$-modules with the corresponding representations
$\rho_U: {\cal L}(R)\rightarrow \End(U)$ and $\rho_W: {\cal L}(R)\rightarrow \End(W)$.
To define the action of the RE algebra
$$
{\cal L}(R)\otimes U\otimes W\rightarrow U\otimes W:\qquad
a\otimes u\otimes w\rightarrow a\triangleright(u\otimes w),
$$
where $a\in{\cal L}(R)$, $u\otimes w\in U\otimes W$, we apply the coproduct
$\Delta(a) = a_{(1)}\otimes a_{(2)}$ (in the Sweedler's notation), then permute
$a_{(2)}$ with the vector $u$ by means of the categorical commutativity
morphism $\mbox{\sf R}_{({\cal L}(R),U)}$:
$$
\mbox{\sf R}_{({\cal L}(R),U)}(a_{(2)}\otimes u) = u_{(3)}\otimes a_{(23)},
$$
and, finally, apply the representations $\rho_U(a_{(1)})$ and $\rho_W(a_{(23)})$
to the corresponding modules:
\be
a\otimes u\otimes w\stackrel{\Delta}{\longrightarrow}a_{(1)}\otimes a_{(2)}\otimes
u\otimes w\stackrel{\mbox{\tiny$\mbox{\sf R}_{{\cal L}(R),U}$}}{\longrightarrow}
a_{(1)}\otimes u_{(3)}\otimes a_{(23)}\otimes w\stackrel{\triangleright}
{\longrightarrow} (a_{(1)}\triangleright u_{(3)})\otimes (a_{(23)}\triangleright w).
\label{gen-ext}
\ee

Below we shall use the following notation for a product of R-matrices:
$$
R^{\pm 1}_{i\rightarrow j} =
\left\{
\begin{array}{lcl}
R^{\pm 1}_iR^{\pm 1}_{i+1}\dots R^{\pm 1}_j&\quad&i<j\\
\rule{0pt}{5mm}
R^{\pm 1}_iR^{\pm 1}_{i-1}\dots R^{\pm 1}_j&\quad&i>j
\end{array}
\right..
$$

Taking the linear combinations
$R_{1\rightarrow p}x_1\otimes\dots x_p$
as the basis vectors of $V^{\otimes p}$
and using (\ref{gen-ext}), (\ref{K-V})
and (\ref{rep-gen}) we get
\begin{eqnarray}
L_1\triangleright R_{1\rightarrow p}x_1\otimes x_2\otimes \dots x_p
&=& (L_1\triangleright R_1x_1)\stackrel{.}{\otimes}L_2\triangleright
(R_{2\rightarrow p}x_2\otimes \dots x_p) =\dots\nonumber\\
&=&\eta^p R^{-1}_{1\rightarrow p}x_1\otimes x_2\otimes \dots \otimes x_p.
\label{act-Vp}
\end{eqnarray}

The chain (\ref{gen-ext}) specialized to $a=L_i^j$, $u=x_k$ leads to an important
consequence. Taking into account  (\ref{K-V}), we find for any $w\in W$
$$
L_1R_1\triangleright (x_1\otimes w) =\eta\, R_1^{-1}x_1\stackrel{.}{\otimes}
(L_2\triangleright w),
$$
or, omitting an arbitrary $w$, we come to the ``permutation rule'' of the operators
$L_i^{\,j}\triangleright $ and basis vectors $x_p$ of the space $V$:
\be
R_1(L_1\triangleright)R_1x_1 = \eta\,x_1(L_2\triangleright).
\label{perm-kx}
\ee
This formula includes the action of $L$ on $V$ and the categorical commutativity
morphism (\ref{K-V}) and gives a simple way of extending the module structure over the RE algebra
to the tensor power $V^{\otimes p}$. It serves us as the key relation for definition of the
BD algebra (\ref{free-qda}) in the next section.

For the dual vector space $V^*$ and its tensor powers the representation
structure is as follows. The representation of the RE algebra in $V^*$ is given by the operators
$$
L_i^j\triangleright y^k = \tilde\eta\,\sum_{s=1}^N y^s (R^2)_{si}^{kj}
$$
or, in a compact matrix form,
\be
L_2\triangleright y_1 = \tilde\eta\,y_1R^2_1,
\label{rep-dual}
\ee
where $\tilde\eta$ is another (nonzero) numerical parameter.

The categorical commutativity morphism (\ref{L-Vd}) and the action (\ref{rep-dual}) leads
to the corresponding operator-vector ``permutation rule'' (analogous to (\ref{perm-kx}))
\be
(L_2\,\triangleright)y_1 = \tilde\eta\, y_1R_1(L_1\triangleright)R_1.
\label{perm-ly}
\ee

There exists a remarkable connection between the set of ${\cal L}(R)$-submodules in
$V^{\otimes p}$ and the R-matrix representation of the Hecke algebra ${\cal H}_p(q)$
in $\End(V^{\otimes p})$.
\begin{proposition}
\label{rea-rep-s}
For any given $p\geqslant 2$ the ${\cal L}(R)$-module $V^{\otimes p}$ is
reducible. The invariant subspaces $V_\lambda\subset
V^{\otimes p}$, $\lambda\vdash p$, are extracted  by the action of projection
operators $P^a_\lambda ={\rho\str{-1.3}}_{\! R}(e^a_\lambda)$, $1\leqslant a\leqslant
d_\lambda $, where
$e^a_\lambda(\sigma)$ is the {\it primitive} idempotent of the Hecke algebra
${\cal H}_p(q)$ corresponding to a standard Young tableau associated with a
partition $\lambda $ (there are $d_\lambda $ of such tableaux in all). Thus, we have
the following expansion:
$$
V^{\otimes p} \cong \bigoplus_{\lambda \vdash p}d_\lambda V_\lambda,\qquad
V_\lambda  \cong {\rm Im}P^a_\lambda,\quad 1\leqslant \forall\,a \leqslant
d_\lambda.
$$
Here the coefficient $d_\lambda $ in the direct sum of vector spaces stands for
the multiplicity of the module $V_\lambda $ in the tensor power $V^{\otimes p}$.
A similar decomposition is true for an ${\cal L}(R)$-module $(V^*)^{\otimes p}$.
\end{proposition}
For more detailed treatment and technical results the reader
is referred to \cite{S}.

To complete the section, we consider another particular module over the RE algebra,
namely, the module $V\otimes V^*$.  As was mentioned above, the corresponding
representation can be treated as the action of the RE algebra ${\cal L}(R)$ on the
generating space $W(L)$ and can be extended to the whole algebra ${\cal L}(R)$
while preserving the algebraic structure of the RE algebra. Due to this reason we
call this representation {\it adjoint}. Such a terminology is also justified by the classical
limit $q\rightarrow 1$ considered below.

So, we consider the ``second copy'' of the space $V\otimes V^*$ and denote its
basis elements as $M_i^j = x_i\otimes y^j$. Thus, the space
$W(M)={\rm span}_{\Bbb K}(M_i^{\,j})$ is isomorphic to the space $W(L)$ generating
the RE algebra and plays the role of the adjoint representation space for the RE algebra.

The commutativity morphism $\mbox{\sf R}_{(W(L),W(M))}$ is given by (\ref{comm-KM})
with the corresponding change of notation for the second factor:
$$
\mbox{\sf R}_{(W(L),W(M))}(L_1\stackrel{.}{\otimes}M_{\overline 2}) =
M_{\overline 2}\stackrel{.}{\otimes}L_1.
$$
Then formulae (\ref{perm-kx}) and (\ref{rep-dual}) allow us to get
the action of the RE algebra ${\cal L}(R)$ on the space $W(M) = {\rm span}_{\Bbb K}\{M_i^j\}$
\be
L_1\triangleright M_{\overline 2} = \eta\,\tilde\eta\,M_{\underline 2}.
\label{action-KM}
\ee
Finally, the adjoint action (\ref{action-KM}) together with the above commutativity
morphism leads to the ``permutation rule'' for the operators $L\,\triangleright $ and
the basis vectors $M$ of the representation space $W(M)$:
\be
(L_1\triangleright) M_{\overline 2} = \eta\,\tilde\eta\, M_{\underline 2}\,(L_1\triangleright).
\label{ad-perm}
\ee
This formula is consistent with the braided bialgebra structure of the RE algebra and
the adjoint action on the space $W(M)$. It gives a way of extending the left
module structure to the whole tensor algebra $TW(M)$.

\section{Braided differential algebras arising from the representation theory of RE algebra}

\label{sec:qd-rea}

In this section we consider the construction of unital associative algebras
${\cal B}({\cal L}(R),{\cal M})$, containing two subalgebras --- a $GL(m)$-type
RE algebra ${\cal L} (R)$ and an ${\cal L} (R)$-algebra ${\cal M}$ which (as a vector
space) is the direct sum of some ${\cal L}(R)$-modules. The subalgebra ${\cal M}$
will be interpreted as a noncommutative function algebra endowed with an action of
 ``exponentiated'' differential operators which form the subalgebra ${\cal L}(R)$.
Due to this reason, we call the algebras ${\cal B}({\cal L}(R),{\cal M})$ the {\it braided
differential algebras} (or BD algebras for short) in what follows. To clarify the reasons
for using such a terminology we consider a classical limit ($q\rightarrow 1$)
of some algebras ${\cal B}({\cal L}(R),{\cal M})$ and suggest the
differential-geometric interpretation of constructions obtained in this way.

In defining the associative algebra structure in ${\cal B}({\cal L}(R),{\cal M})$ a
decisive role belongs to the permutation rule of elements of ${\cal L}(R)$ and
${\cal M}$. This should be an analog of the classical Leibnitz rule, since it
embraces the action of a differential operator on
a function and their mutual permutation (see (\ref{perm-kx}), (\ref{perm-ly}) and
(\ref{ad-perm}) as examples). We shall refer to this rule as the {\it operator-function
permutation} (OFP) rule.

We impose two natural requirements on the OFP rule. First, it should respect the
algebraic structures of ${\cal L}(R)$ and ${\cal M}$ as subalgebras of the BD algebra.
This means, that the subalgebra ${\cal M}$ is an ${\cal L}(R)$-module and the
action of the RE algebra  ${\cal L}(R)$ is compatible with the multiplication in
$\cal M$ (that is $\cal M$ is an ${\cal L} (R)$-algebra).

Second, the OFP relation must be compatible with possible additional
symmetries of ${\cal L}(R)$ and $\cal M$. As an example of such a symmetry we
can point out the coadjoint comodule structure of ${\cal L}(R)$ over the RTT algebra
(see (\ref{rtt-rea-coac})). The ``quantum function" algebra $\cal M$ can also bear
coadjoint or (co)vector comodule structure over the RTT algebra\footnote{The last case
was considered in the paper \cite{IP} devoted to the Heisenberg double over
the quantum group.}.

It turns out that the first requirement restricts considerably possible forms of the OFP relation.
Besides, the RE algebra representation theory and the structure
of the Schur-Weyl category allows one to find all possible OFP rules up to a
renormalization isomorphism.

Below we give several important examples of BD algebras. In the next section
we use them in order to construct a BD algebra involving general quantum {\em matrix} algebra
${\cal M}$ and the RE algebra acting on ${\cal M}$ by quantum right-invariant differential operators.
\medskip

\nonumber {\bf Example 1.} Let an $N$-dimensional vector space $V$ be a left
${\cal L}(R)$-module with the action (\ref{rep-gen}) of the ${\cal L}(R)$ generators on a given
basis set $\{x_i\}_{1\leqslant i\leqslant N}$ of the space $V$. Consider a unital
associative $\Bbb K$-algebra ${\cal X}(V)$ freely generated by elements $x_i$:
$$
{\cal X}(V) = {\Bbb K}\langle x_1,x_2,\dots ,x_N\rangle
$$
and its $p$-th order homogeneous component ${\cal X}^p(V) \simeq V^{\otimes p}$. The
${\cal L}(R)$-module structure is introduced by an
analog of the relation (\ref{perm-kx}).
This formula is the key point for constructing the BD algebra
${\cal B}({\cal L}(R),{\cal X}(V))$ --- it leads to the OFP relation we need.
\begin{definition}{\rm
Let ${\cal X}(V)={\Bbb K}\langle x_i\rangle_{1\leqslant i\leqslant N}$ be an algebra of
noncommutative polynomials freely generated by elements $x_i$, ${\cal L}(R)$ be
the RE algebra generated by $N^2$ elements $L_i^{\,j}$ subject to multiplication rules
(\ref{REA}) with a $GL(m)$-type R-matrix. Then the {\it free braided differential algebra}
is the unital associative algebra ${\cal B}({\cal L} (R),{\cal X}(V))$ generated
by $\{x_i\}$ and $\{L_i^{\, j}\}$ subject to the additional permutation  rule
\be
 R_1L_1R_1x_1 = \eta\,x_1L_2,\qquad \eta\in {\Bbb K}\setminus 0.
\label{free-qda}
\ee }
\end{definition}

In order to provide the subalgebra ${\cal X}(V)$ with the structure of a module over the RE algebra,
we should only define an action of $L$ generators on the unit element
$1_\sscr{$\,\cal B$}$. Since this action should realize a one-dimensional
representation of the RE algebra, we naturally set
\be
L\triangleright 1_\sscr{$\,\cal B$} = \varepsilon(L)\, 1_\sscr{$\,\cal B$}.
\label{unit-ac}
\ee

Then OFP relation (\ref{free-qda}) together with (\ref{unit-ac}) allows us to get the
action of $L$ on any homogeneous monomial in $x_i$: we should move the element $L$
to the most right position and then apply (\ref{unit-ac}). For example, for a $p$-th order
homogeneous monomial we find
\begin{eqnarray*}
(R_{p\rightarrow 1}&\hspace*{-3.9mm}L_1&\hspace*{-3.9mm}R_{1\rightarrow p})
\triangleright( x_1x_2\dots x_p)\equiv (R_{p\rightarrow 1}L_1R_{1\rightarrow p}\,x_1
x_2\dots x_p)\triangleright 1_{\sscr{$\,\cal B$}}\\
&\hspace*{-3.9mm}=\hspace*{-3.9mm}&( R_{p\rightarrow 2}(R_1L_1R_1x_1)
R_{2\rightarrow p}\,x_2\dots x_p)\triangleright 1_{\sscr{$\,\cal B$}}\stackrel{(\ref{free-qda})}{=}
\eta x_1(R_{p\rightarrow 3}(R_2L_2R_2x_2) R_{3\rightarrow p}\,x_3\dots x_p)
\triangleright 1_{\sscr{$\,\cal B$}}\\
&\hspace*{-3.9mm}=\hspace*{-3.9mm}&\dots = \eta^p x_1x_2\dots x_p(L_{p+1}\triangleright
1_{\sscr{$\,\cal B$}}) \stackrel{(\ref{unit-ac})}{=} \eta^px_1x_2\dots x_p\,I_{p+1}.
\end{eqnarray*}
Clearly, this is the same action as (\ref{act-Vp}) in full agreement with the isomorphism
${\cal X}^p(V)\simeq V^{\otimes p}$.

It is evident, that the free BD algebra ${\cal B}({\cal L}(R),{\cal X}(V))$ contains all the ${\cal L}(R)$-modules
$V_\lambda $, $\lambda\vdash p\geqslant 1$. Any such a module is a subspace
of the corresponding homogeneous component ${\cal X}^p(V)$:
$$
V_\lambda\cong {\rm Im}\Big({\rho\str{-1.3}}_{\! R}(e_\lambda^a)\Big)\subset
{\cal X}^p(V), \quad \lambda\vdash p,\quad 1\leqslant a\leqslant d_\lambda.
$$
with the multiplicity $d_\lambda$ (see Proposition \ref{rea-rep-s}).

We can decrease the size of the free BD algebra by passing to a quotient
$$
{\cal X}_J (V)= {\cal X}(V)/\langle J\rangle, \quad J\subset {\cal X}(V).
$$
Recall, that $\langle J \rangle $ stands for the two-sided ideal, generated by a subset
$J$. Assuming  the ideal $\langle J\rangle $ to be invariant w.r.t. the action of
${\cal L}(R)$  we can define its action on the quotient  ${\cal X}_J(V)$.

A systematic way to get a set of relations on $x_i$ with the desired properties
consists in choosing $J$  to be equal to the image of a {\it central} idempotent
$e_\lambda(\sigma)\in {\cal H}_p(q)$ for some $p\geqslant 2$:
$$
J_\lambda= {\rm Im}\Big({\rho\str{-1.3}}_{\! R}(e_\lambda)\Big)\subset {\cal X}
^p(V), \quad
\lambda\vdash p.
$$
Basing on the properties of idempotents $e_\lambda $ one can show that at the
canonical projection $\pi_\lambda: {\cal X}(V)\rightarrow {\cal X}_{J_\lambda}(V)$
all the ${\cal L}(R)$-submodules $V_\mu\in {\cal X}(V)$ corresponding to partitions
$\mu\supset \lambda $ are mapped to zero:
$$
\pi_\lambda(V_\mu) = 0,\quad\forall \mu\supset \lambda.
$$

For example, if we want to impose  {\it quadratic} relations on the generators
$x_i$ we have only two possibilities: to annihilate the $q$-antisymmetric
component
\be
J_{(1^2)}\subset {\cal X}^2(V):\qquad J_{(1^2)} = {\rm Im}\Bigr((q-R)\Bigl)
\label{J-anti}
\ee
or $q$-symmetric component
\be
J_{(2)}\subset
{\cal X}^2(V):\qquad J_{(2)} = {\rm Im}\Bigr((q^{-1} +R)\Bigl).
\label{J-symm}
\ee
The choice (\ref{J-anti}) gives rise to a BD algebra ${\cal B}({\cal L}(R),{\cal X}_{s}(V))$
of the RE algebra ${\cal L}(R)$ over  the ``quantum plane'' ${\cal X}_{s}(V)$ \cite{FRT}
\be
\begin{array}{l}
R_1x_1x_2 - q x_1x_2 =0\\
\rule{0pt}{5mm}
R_1L_1R_1L_1 - L_1R_1L_1R_1 = 0\\
\rule{0pt}{5mm}
R_1L_1R_1x_1 = \eta\, x_1L_2.
\end{array}
\label{rea-plane}
\ee
The action ${\cal L}(R)\triangleright {\cal X}_s(V)$ is induced by (\ref{unit-ac})
together with the third relation in system (\ref{rea-plane}). The BD algebra (\ref{rea-plane})
contains only the ${\cal L}(R)$-submodules isomorphic to $V_{(p)}$, where $(p)$ is a
single row partition of an integer $p\geqslant 1$.

The BD algebra ${\cal B}({\cal L}(R),{\cal X}_{s}(V))$ is covariant with respect to the left
coaction of the RTT bialgebra
$$
\delta_\ell(L_i^{\,j}) = \sum_{k,p=1}^N T_i^{\,k}S(T_p^{\,j})\otimes L_k^{\,p},\qquad
\delta_\ell(x_i) = \sum_{k=1}^NT_i^{\,k}\otimes x_k,
$$
$S$ being the antipodal map in the RTT algebra. In practical calculations it is convenient
to use a loose notation $\delta_\ell(L_1) = T_1L_1S(T_1)$ and $\delta_\ell(x_1) = T_1x_1$ and treat
$T_a$ to be commutative with $L_b$ and $x_c$ if $a\not= b$ and $a\not=c$.

Let us verify the covariance of the OFP rule in the system (\ref{rea-plane}). We get
\begin{eqnarray*}
\delta_\ell(R_1L_1R_1x_1) &=& R_1T_1L_1\underline{S(T_1)R_1T_1}x_1
=R_1T_1L_1T_2R_1S(T_2)x_1 = \underline{R_1T_1T_2}L_1R_1x_1S(T_2) \\
&=&T_1T_2\underline{R_1L_1R_1x_1}S(T_2) = \eta\, T_1x_1
T_2L_2S(T_2) = \delta_\ell(\eta\, x_1L_2).
\end{eqnarray*}

Assuming a given R-matrix of $GL(m)$-type to be a deformation of the usual flip $P$
(then $m=N=\dim_{\Bbb K} V$),
consider the classical limit $q\rightarrow 1$ of the BD algebra ${\cal B}({\cal L}
(R),{\cal X}_{s}(V))$ (\ref{rea-plane}). For this purpose we pass to a different set
$\{K_i^{\,j}\}_{1\leqslant i,j,\leqslant N}$ of the RE algebra generators:
\be
L = I-(q-q^{-1})K,\qquad K=\|K_i^{\,j}\|.
\label{shift}
\ee
Taking into account the Hecke condition (\ref{Hec}), we rewrite the defining
relations (\ref{REA}) in terms of the new generators
\be
R_1K_1R_1K_1 - K_1R_1K_1R_1 = R_1K_1 - K_1R_1.
\label{mrea}
\ee
The bialgebra structure now reads
\be
\Delta(K) = 1\otimes K + K\otimes 1 - (q-q^{-1}) K\otimes K,\qquad
\varepsilon(K) = 0 .
\label{new-copr}
\ee

Then, according to the first line of (\ref{rea-plane}), the generators $x_i$ of the
subalgebra ${\cal X}_s(V)$ turn into commutative elements
\be
x_2x_1 - x_1x_2 = 0.
\ee
So, at $q\rightarrow 1$ we have ${\cal X}_s(V)={\Bbb K}[V^*]$.

The multiplication rules (\ref{mrea}) turns into defining relations of the universal
enveloping algebra $U(gl(m))$
\be
\kappa_1\kappa_2 - \kappa_2\kappa_1 = \kappa_1P_{12} - P_{12}\kappa_1,
\label{gl-lim}
\ee
where the matrix $\kappa = \|\kappa_i^{\,j}\|$ is the limit of generating matrix $K$ at
$q\rightarrow 1$.

In order to get the limit of the RE algebra action (the third relation in (\ref{rea-plane})) we additionally
suppose the following behavior of the parameter $\eta = 1-(q-q^{-1})\eta_0+o(q^2-1)$. Under
this assumption the OFP relation in (\ref{rea-plane}) gives rise to
\be
\kappa_2 x_1 - x_1\kappa_2= \eta_0x_1 +P_{12}x_1.
\label{ell-act}
\ee
Together with the commutation relations (\ref{gl-lim}) this formula allows us to interpret the
generators $\kappa_i^{\,j}$ as the following vector fields on the ${\Bbb K}[V^*]$:
\be
\kappa_i^{\,j} = x_i{\partial\hspace{0.8pt}}^j_x+\eta_0\delta_i^j\,(x\cdot\partial_x),
\label{v-field}
\ee
where we denote
$$
{\partial\hspace{0.8pt}}^k_x=
\frac{\partial}{\partial x_k},\qquad (x\cdot \partial_x) =\sum_{k=1}^m
x_k{\partial\hspace{0.8pt}}^k_x.
$$
If $V$ is the left fundamental vector $GL(m)$-module
$$
x_i\mapsto x_j M^j_{\,i},\quad M = \|M^j_{\,i}\|\in GL(m)
$$
then the fields $\kappa_i^{\,j}$ in (\ref{v-field}) are invariant with respect
to the $GL(m)$ action on the right side.
\begin{remark}
{\rm
In considering the classical limit $q\rightarrow 1$ it is convenient to parameterize
$q=e^{{\tau\over 2}}$ and treat the classical limit as $\tau\rightarrow 0$.
In this limit the shift formula (\ref{shift}) turns into $L = I -\tau\kappa +o(\tau^2)$.
Together with the group-like coproduct (\ref{copr}) for $L$ generators and their
``Weyl-type'' commutation with the generators $x_i$ (the third relation in system
(\ref{rea-plane})) it allows us to interpret the generators $L_i^{\,j}$ as
exponentiated quantized differential operators $\kappa$ (\ref{v-field}).

Note, that $x\cdot \partial = {\rm Tr}\,\kappa $ is a {\it central} element of the Lie
algebra $gl(m)$. Therefore, on adding to $\kappa_i^{\,j}$ (\ref{v-field}) a term
proportional to this central element, we can specialize the parameter $\eta_0$ in
(\ref{ell-act}) to any given value (for example, we can get $\eta_0=0$). Such an
operation changes the multiplicative parameter $\eta $ in the OFP relation of
(\ref{rea-plane}). This is another evidence of exponential-like dependence of
quantum differential operators $L$ on classical differential operators $\kappa $:
a linear shift of $\kappa $ leads to a multiplicative renormalization of $L$.
}
\end{remark}

Consider now the choice (\ref{J-symm}) for the permutation rules on $x_i$.
We come to the BD algebra ${\cal B}({\cal L}(R), {\cal X}_a(V))$ with the following
relations on the generators
\be
\begin{array}{l}
R_1x_1x_2 +q^{-1} x_1x_2 =0\\
\rule{0pt}{5mm}
R_1L_1R_1L_1 - L_1R_1L_1R_1 = 0\\
\rule{0pt}{5mm}
R_1L_1R_1x_1 = \eta\, x_1L_2.
\end{array}
\label{ext-plane}
\ee
This algebra has the same RTT-comodule property as the algebra
${\cal B}({\cal L}(R), {\cal X}_s(V))$ considered above. But contrary to the BD algebra
(\ref{rea-plane}), now we have only finite number of the ${\cal L}(R)$-submodules in
${\cal B}({\cal L}(R), {\cal X}_a(V))$. Namely, the vector space ${\cal X}_a(V)$ is isomorphic
to the direct sum of the modules $V_{(1^p)}$, $1\leqslant p\leqslant m$.

Again, assuming  a given R-matrix $R$ to be a deformation of the flip $P$ and applying the shift (\ref{shift}) we get
at the classical limit $q\rightarrow 1$ the following system of relations
$$
\begin{array}{l}
x_1x_2 +x_2x_1 =0\\
\rule{0pt}{5mm}
\kappa_1\kappa_2 - \kappa_2\kappa_1 = \kappa_1P_{12} - P_{12}\kappa_1\\
\rule{0pt}{5mm}
\kappa_2 x_1 - x_1\kappa_2= \eta_0x_1 +P_{12}x_1.
\end{array}
$$

The classical algebra defined by the above relations on the generators has a transparent
geometrical interpretation.
The elements $x_i$ generate an external subalgebra and are treated
as one-forms --- the differentials of coordinate functions of ${\Bbb K}[V^*]$
$$
x_i = dy_i,\quad 1\leqslant i\leqslant m,
$$
while $\kappa_i^j$ is of the form
$$
\kappa_i^j = \frak{L}_i^j + \eta_0\,\delta_i^j\sum_{k=1}^m\frak{L}_k^k,
$$
where $\frak{L}_i^j$ is the {\it Lie derivative} along the vector field
$y_i{\partial\hspace{0.8pt}}_y^j$.
\medskip

{\bf Example 2.}
We can start from a more interesting  module $V\otimes V^*$ (called {\em adjoint}) with
the linear basis $M_i^{\, j} = x_i\otimes y^j$. The action of the RE algebra  is given by (\ref{action-KM}).
Formula (\ref{ad-perm}) provides a recipe for extending the module structure on
tensor powers of the adjoint module $V\otimes V^*$. In analogy with constructions of
Example 1, we consider a unital associative algebra ${\cal M}$, generated by $N^2$
free elements $M_i^{\,j}$ and define the algebra ${\cal B}({\cal L}(R), {\cal M})$ by
imposing the following multiplication rules of the free generators $M_i^{\,j}$ and RE algebra
generators $L_i^{\,j}$
\be
L_1M_{\overline 2} = M_{\underline 2}L_1.
\label{KM-d}
\ee
This formula stems from the relation (\ref{ad-perm}) which is defined by the RE algebra
representation theory. Then the subalgebra ${\cal M}\subset {\cal B}({\cal L}(R),{\cal M})$
can be given the structure of a module over the RE algebra by relation (\ref{unit-ac}).
Note, that the requirement (\ref{unit-ac}) fixes the constant $\tilde\eta$ as $\eta\,\tilde\eta = 1$.

We can restrict the algebra ${\cal B}({\cal L}(R),{\cal M})$ by setting some
relations on the generators $M_i^{\,j}$ which are consistent with the OFP relation
(\ref{KM-d}). We consider the case of {\it quadratic} relations. In \cite{GPS1}
a pair of orthogonal projectors ${\cal A}_q, {\cal S}_q: {\cal M}^{(2)}\rightarrow
{\cal M}^{(2)}$ was constructed. Here ${\cal M}^{(2)}$ is the subspace of ${\cal M}$
spanned by the quadratic monomials in generators $M_i^{\,j}$. The projectors
${\cal A}_q$ and ${\cal S}_q$  have the natural interpretation as a
$q$-antisymmetrizer and a $q$-symmetrizer on the space ${\cal M}^{(2)}$.
The images of these operators are invariant subspaces with respect
to the RE algebra action. So, a consistent quadratic relation on the free generators $M_i^{\,j}$
can be chosen as ${\rm Im}\,{\cal A}_q = 0$ or ${\rm Im}\,{\cal S}_q=0$.

Consider the first case. It can be shown that the requirement ${\rm Im}\,{\cal A}_q
= 0$ is equivalent to the RE algebra type relations on the generators $M_i^{\,j}$. So, we come
to the BD algebra  defined by the following relations on generators
\be
\begin{array}{l}
R_1M_1R_1M_1 - M_1R_1M_1R_1=0\\
\rule{0pt}{5mm}
R_1L_1R_1L_1 - L_1R_1L_1R_1 = 0\\
\rule{0pt}{5mm}
R_1L_1R_1M_1 =  M_1R_1L_1R_1.
\end{array}
\label{2-REA-d}
\ee
We denote this BD algebra ${\cal B}_{\rm ad}({\cal L}(R),{\cal M}(R))$.
Both the RE algebras ${\cal L}(R)$ and ${\cal M}(R)$ are subalgebras of
${\cal B}_{\rm ad}({\cal L}(R),{\cal M}(R))$, the subalgebra ${\cal M}(R)$ is endowed with
a ${\cal L}(R)$-module structure by means of (\ref{unit-ac}) and by the third relation of
the above system. Besides, the algebra ${\cal B}_{\rm ad}({\cal L}(R), {\cal M}(R))$ has the
left coadjoint comodule structure over the RTT algebra (\ref{RTT}).

The algebraic properties of this BD algebra will be considered in more detail in the next
section. Here we only point out that the R-traces ${\rm Tr\str{-1.3}}_{\!R}M^k$,
$k\geqslant 0$, are central in the whole BD algebra ${\cal B}_{\rm ad}({\cal L}(R),{\cal M}(R))$
(not only in the RE subalgebra ${\cal M}(R)$) and therefore are invariant under
the action of the subalgebra ${\cal L}(R)$. This means that this action can be restricted
to quotients of ${\cal M}(R)$ over ideals generated by relations
${\rm Tr\str{-1.3}}_{\!R}M^k = c_k$, $1\leqslant k\leqslant m$, where $c_k$ are
fixed constants. Recall, that in \cite{GS1, GS4} such like quotients were interpreted as
quantum (braided) analogs of  $GL(m)$ orbits  in ${gl^*(m)}$. Therefore, the
subalgebra ${\cal L}(R)$ in the BD algebra ${\cal B}_{\rm ad}({\cal L}(R),{\cal M}(R))$ can be treated
as the quantized algebra of differential operators generated by the vector fields
tangential to the mentioned orbits.

To justify this interpretation we consider the classical limit $q\rightarrow 1$
of the BD algebra (\ref{2-REA-d}) by assuming  $R$ to be a deformation of the flip $P$.

Making the shift (\ref{shift}) for $L$ and passing to the limit $q\rightarrow 1$ in the BD algebra
(\ref{2-REA-d}) we come to the following permutation rules for the generators
$m_i^{\,j} = \lim_{q\rightarrow 1} M_i^{\,j}$ and $\kappa_i^{\,j} = \lim_{q\rightarrow 1} K_i^{\,j}$:
$$
\begin{array}{l}
m_1m_2 -  m_2m_1 = 0\\
\rule{0pt}{5mm}
\kappa_1\kappa_2 -\kappa_2\kappa_1 = \kappa_1P_{12} - P_{12}\kappa_1\\
\rule{0pt}{5mm}
\kappa_2m_1-m_1\kappa_2 = P_{12}m_1-m_1P_{12}.
\end{array}
$$
The two last lines in this system
of permutation rules show that $\kappa_i^{\,j}$ are coadjoint vector fields on
the space $gl^*(m)$:
\be
\kappa_i^{\,j} = m_i^{\,s}\frac{\partial}{\partial m_{j}^{\,s}} -
m_s^{\,j}\frac{\partial}{\partial m_{s}^{\,i}},
\label{coa-vf}
\ee
where the summation over the index $s$ is understood. As is well known,
the vector fields (\ref{coa-vf}) are tangent to the $GL(m)$ orbits in the linear space $gl^*(m)$.

\section{The braided differential algebra over QM algebra}
\label{sec:5}

In the preceding section we gave an example (\ref{2-REA-d}) of a BD algebra over a
quantum matrix algebra. In the classical limit the RE algebra generators turned
into the adjoint vector fields (\ref{coa-vf}). In the example (\ref{2-REA-d})
the algebra of "quantum functions" was taken to be the second copy
of RE algebra.

Now we are going to define the braided differential algebra ${\cal B}_r({\cal L}(R),{\cal M})$
of the $GL(m)$-type RE algebra ${\cal L}(R)$ over an arbitrary quantum matrix algebra
${\cal M}(R,F)$. Therefore, we should supply the BD algebra with the OFP relation similar
to the third relation in (\ref{rea-plane}). The construction is presented in the following definition.

\begin{definition}
\label{def:bda-gen}
{\rm
Let ${\cal L}(R)$ be the RE algebra associated with a $GL(m)$-type R-matrix
$R$ and ${\cal M}(R,F)$ be the QM algebra, associated with a compatible pair
of R-matrices $\{R,F\}$ (see section \ref{sec:hr-qma}). Define a unital
associative algebra ${\cal B}_r({\cal L}(R),{\cal M})$ over the field $\Bbb K$
generated by the elements $L_i^{\,j}$ of the RE algebra and by elements $M_i^{\,j}$
of the QM algebra subject to the following system of relations
\begin{eqnarray}
&&
R_1M_1M_{\overline 2} - M_1M_{\overline 2}R_1 = 0
\nonumber\\
&&
R_1L_1R_1L_1 - L_1R_1L_1R_1 = 0 \nonumber\\
&&
R_1L_1R_1M_1 = \eta M_1L_{\overline 2},
\label{BDA-gen-left}
\end{eqnarray}
where the ``matrix copies" $M_{\overline 2}$ and $L_{\overline 2}$
are produced by the R-matrix $F$ in accordance with (\ref{copy}):
$$
L_{\overline 2} = F_1L_1F_1^{-1},\qquad M_{\overline 2} = F_1M_1F_1^{-1}.
$$
The nonzero number $\eta$ is a parameter of the algebra.

We introduce an action of the RE algebra generators on the unit element
$1_\sscr{$\,\cal B$}$ by the rule
\be
a \triangleright 1_\sscr{$\,\cal B$} = \varepsilon(a) 1_\sscr{$\,\cal B$},
\quad \forall a\in {\cal L}(R),
\label{rea-action}
\ee
where $\varepsilon$ is the counit map of the braided bialgebra ${\cal L}(R)$:
$\varepsilon(L) = I$ (see (\ref{coun})).
}
\end{definition}

Note, that the Heisenberg double, considered in \cite{IP} corresponds to the pair $\{R,P\}$
of the compatible R-matrices, where $R$ is a $GL(m)$-type R-matrix. In this case the QM
algebra ${\cal M}(R,P)$ turns into the Hopf algebra of quantum functions on $GL(m)$.

\begin{remark}{\rm
Let us shortly explain how one can get the OFP relation (\ref{BDA-gen-left}) in the above
definition. As we mentioned in section \ref{sec:rea-rep}, the RE algebra generators can be
treated as basis elements of the space $V\otimes V^*$ of the Schur-Weyl category {\sf SW(V)}
(the category of finite dimensional modules over the RE algebra): $L_i^{\,j} = x_i\otimes y^j$,
$\{x_i\}_{ 1\leqslant i\leqslant N}$ and $\{y^j\}_{ 1\leqslant j\leqslant N}$ being the respective
basis (\ref{duba}) of $V$ and $V^*$.  In order to get the RE algebra action on generators $M_i^{\,j}$
which would be not the adjoint one (as in (\ref{2-REA-d})) but rather similar to (\ref{rea-plane})
we proceed as follows.

Let us enlarge the class of objects of the category {\sf SW(V)} by another pair of $N$ dimensional
vector spaces
$$
U = {\rm span}_{\Bbb K}(t_i)_{ 1\leqslant i\leqslant N},\qquad
U^* = {\rm span}_{\Bbb K}(z^i)_{1\leqslant i\leqslant N}
$$
dual to each other. Then the matrix elements $M_i^{\,j}$ are taken
to be the basis elements of $V\otimes U^*$: $M_i^{\,j} = x_i\otimes z^j$.
In order to get the OFP relation among $L_i^j$ and $M_r^s$ we have
to take into account the (known) OFP relation among $L$ and $x$ and the
categorical permutation morphism of $L$ (treated as $V\otimes V^*$) and $U^*$
in accordance with our general recipe described in (\ref{gen-ext}).

The commutativity morphism $\mbox{\sf F}:V\otimes U^*\rightarrow U^*\otimes V$
is defined via the operator $F$:
$$
 \mbox{\sf F}(x_i\otimes z^j) = z^k\otimes x_sF_{ki}^{\, sj}.
$$
This gives rise to $F$ appearing in (\ref{BDA-gen-left}) in the formula for $L_{\overline 2}$.

Note also, that there exists another choice of the ``function algebra''. Namely,
we could take as the basis of function algebra the elements
$M_i^{\,j} = t_i\otimes y^j\in U\otimes V^*$. Such a choice would give
rise to another form of the OFP relation (compare with (\ref{BDA-gen-left}))
\be
L_{\overline 2}M_1 = \tilde\eta\, M_1R_1L_1R_1.
\label{BDA-gen-right}
\ee

Actually, we can get the corresponding BD algebra starting from (\ref{BDA-gen-left}).
If we introduce the matrix $\hat L = M^{-1}LM$ and the new R-matrix $\hat R = F^{-1}R^{-1}F$,
then the system of relations (\ref{BDA-gen-left}) leads to
\begin{eqnarray*}
&&
\hat R_1M_{\underline 2} M_{\underline 1} - M_{\underline 2} M_{\underline 1}\hat R_1= 0
\\
&&
\hat R_1\hat L_1\hat R_1\hat L_1 - \hat L_1\hat R_1\hat L_1\hat R_1 = 0 \\
&&
\hat L_{\underline 2} M_1  =  \eta M_1\hat R_1\hat L_1\hat R_1,
\end{eqnarray*}
where $M_{\underline 2} = F^{-1}M_1F$. The OFP relation standing in the third line of the above
system coincides with (\ref{BDA-gen-right}) (up to the nonessential change $F\rightarrow F^{-1})$.

Strictly speaking, the BD algebra generated by $\hat L$ and $M$ is not a subalgebra of the
algebra (\ref{BDA-gen-left}) since we have to use the inverse matrix $M^{-1}$ in passing from
$L$ to $\hat L$. As follows from the Cayley-Hamilton identity (\ref{qma-ch}) this requires some
extension of the initial algebra, namely, we have to demand the invertibility of the element $a_m(M)$
(see \cite{IP} for more detail).
}
\end{remark}

From the viewpoint of the representation theory of the RE algebra, the BD algebra introduced in
definition \ref{def:bda-gen} consists of the direct sum of modules over the RE algebra isomorphic
to those of the BD algebra (\ref{rea-plane}). To be more precise, the following proposition
holds true.
\begin{proposition}
Relation (\ref{rea-action}) allows us to define the ${\cal L}(R)$-module structure on the
subalgebra ${\cal M}(R,F)$ of the BD algebra ${\cal B}_r({\cal L}(R),{\cal M})$ introduced
in definition \ref{def:bda-gen}. The action of the RE algebra generators $L_i^{\, j}$ on the basis
vectors of the $p$-th order homogeneous component ${\cal M}^p(R,F)$ reads
\be
L_1\triangleright R_{(1\rightarrow p)}M_{\overline 1}M_{\overline 2}\dots
M_{\overline p} = \eta^p R^{-1}_{(1\rightarrow p)} M_{\overline 1}M_{\overline 2}\dots
M_{\overline p},
\label{p-th-ord}
\ee
where $M_{\overline{k}} = F_{k-1}M_{\overline {k-1}}F^{-1}_{k-1}$.
\end{proposition}
\smallskip

\noindent
{\bf Proof.}
The proof consists in direct calculations. First of all, using the compatibility
conditions (\ref{cc}) we can transform relation (\ref{BDA-gen-left}) to
$$
R_k L_{\overline k}R_kM_{\overline k} = \eta\, M_{\overline k}L_{\overline{k+1}},
\quad \forall k\geqslant 1,
$$
where the copies $L_{\overline k}$ and $M_{\overline k}$ are defined with the
help of the R-matrix $F$ in accordance with (\ref{copy}). Then we get{}
$$
L_1\triangleright (R_{(1\rightarrow p)}M_{\overline 1}\dots M_{\overline p}) = \eta^p\,
R^{-1}_{(1\rightarrow p)}M_{\overline 1}\dots M_{\overline p}(L_{\overline{p+1}}
\triangleright 1_{\cal B}).
$$
Since $L_{\overline{p+1}} \triangleright 1_{\cal B} = \varepsilon(L_{\overline{p+1}})1_{\cal B}
= I_{12\dots p+1}1_{\cal B}$, the result (\ref{p-th-ord}) follows. This should be
compared with (\ref{act-Vp}).

In a similar manner one can prove that the RE algebra action respects the algebraic
structure of the QM algebra ${\cal M}(R,F)$, that is
$$
a\triangleright(R_kM_{\overline k}M_{\overline{k+1}} - M_{\overline k}
M_{\overline{k+1}} R_k) = 0,\quad \forall a\in{\cal L}(R),\; \forall\,k\geqslant 1.
\eqno{\rule{6.5pt}{6.5pt}}
$$
\medskip

Now, we consider  the case of the BD algebra over the RE algebra in more detail. This means
that we put $F=R$. Note, that we do not come to the algebra (\ref{2-REA-d}) since the OFP
relation (\ref{BDA-gen-left}) takes the form
\be
R_1L_1R_1M_1 = \eta M_1R_1L_1R_1^{-1}
\label{right-inv-bda}
\ee
which differs from the third relation of the BD algebra (\ref{2-REA-d}) by the inverse $R$
in the last place.
As a consequence, the traces ${\Tr\str{-1.3}}_{\!R}M^k$ are not central in the
BD algebra (\ref{right-inv-bda}), and
the action of the braided differential operators from the subalgebra ${\cal L}(R)$
does not preserve
the quantum orbits which are quotients of the RE algebra ${\cal M}(R)$ over the ideals
generated by
conditions on these traces \cite{GS1}. It is not a surprise, since as can be easily
seen from the
classical limit $q\rightarrow 1$, the relation (\ref{right-inv-bda}) defines the right-invariant vector
fields on the $gl^*(m)$:
$$
L=I-(q-q^{-1})K,\qquad K_i^{\,j}\;\stackrel{q\rightarrow 1}{\longrightarrow}\;
m_i^{\,a}\frac{\partial}{\partial m_j^{\,a}}.
$$
Here we neglect the possible central term proportional to $\eta_0$ (see
(\ref{v-field})).

The BD algebra (\ref{2-REA-d}) consisting of the quantum differential operators
generated by coadjoint vector fields can be subtracted as a
subalgebra of the algebra ${\cal B}_r({\cal L}(R),{\cal M}(R))$ with the OFP relation
(\ref{right-inv-bda}). (More precisely, we extend this algebra by $L^{-1}$ and $M^{-1}$.)

Let us introduce the matrices
\be
Q = LM^{-1}L^{-1}M, \qquad N = M^{-1}Q.
\label{def:Q-N}
\ee
Here we have to impose the invertibility condition on the elements $a_m(L)$ and $a_m(M)$
(the polynomial $a_m$ is defined in (\ref{elem-f})) and extend our algebra by the elements
$a_m^{-1}(L)$ and $a_m^{-1}(M)$. Then, the Cayley-Hamilton identities (\ref{CH-id}) for $L$
and $M$ guarantee the invertibility of the matrices involved. The following proposition is a
direct consequence of the multiplication rule (\ref{right-inv-bda}).
\begin{proposition}
The matrix elements of $Q$ and $M$ satisfy the following multiplication rules
\be
\begin{array}{l}
R_1M_1R_1M_1 - M_1R_1M_1R_1 = 0\\
\rule{0pt}{5mm}
R_1Q_1R_1Q_1 - Q_1R_1Q_1R_1 = 0\\
\rule{0pt}{5mm}
R_1Q_1R_1M_1 - M_1R_1Q_1R_1 = 0.
\end{array}
\label{Q-M}
\ee
For the pair $M$ and $N$ we have
\be
\begin{array}{l}
R_1M_1R_1M_1 - M_1R_1M_1R_1 = 0\\
\rule{0pt}{5mm}
R_1N_1R_1N_1 - N_1R_1N_1R_1 = 0\\
\rule{0pt}{5mm}
R_1^{-1}N_1R_1M_1 - M_1R_1^{-1}N_1R_1 = 0.
\end{array}
\label{N-M}
\ee
\end{proposition}

It is clear from the above proposition that the $Q$ and $M$ generate a subalgebra in the BD
algebra ${\cal B}_{r}({\cal L},{\cal M})$ (\ref{BDA-gen-left}). Moreover, we can calculate the
action of the generators $Q_i^{\,j}$ on the unit element $1_\sscr{$\,\cal B$}$ and turn the above
subalgebra into a new BD algebra.

\begin{proposition}
Given the action (\ref{unit-ac}) of the generators $L$, for the elements
$Q$ defined in (\ref{def:Q-N}) one gets:
\be
Q_1\triangleright 1_{\cal B} = \xi\, I_11_{\cal B},\qquad \xi =
\eta^{-1}q^{2m}.
\label{Q-unit}
\ee
\end{proposition}

We call the BD algebra, generated by $Q$ and $M$ subject to the system of relations
(\ref{Q-M}) and the action (\ref{Q-unit}) the adjoint BD algebra and denote it as
${\cal B}_{\rm ad}({\cal Q},{\cal M})$. Note, that OFP relation in the adjoint BD algebra
does not depend on the parameter $\eta $ entering the OFP relation in the algebra
(\ref{BDA-gen-left}). This parameter appears only in the action of the adjoint generators
$Q_i^{\,j}$ on the unit element in (\ref{Q-unit}).

It is evident that the BD algebra
${\cal B}_{\rm ad}({\cal Q},{\cal M})$ defined by relations (\ref{Q-M}) differs from
the algebra ${\cal B}_{\rm ad}({\cal L}(R),{\cal M}(R))$ of example (\ref{2-REA-d}) only by
notation and general normalization (the parameter $\xi $) of the ${\cal L}(R)$-action.

Based on this result, we can reveal the structure of the ${\cal B}_{\rm ad}({\cal Q},{\cal M})$
as a module over the RE algebra.

\begin{proposition}
Given the adjoint BD algebra ${\cal B}_{\rm ad}({\cal Q},{\cal M})$ defined
by relations (\ref{Q-M}) and (\ref{Q-unit}), the subalgebra ${\cal M}(R)$
generated by $M_i^{\,j}$ is endowed with the ${\cal Q}(R)$-module structure
with the following action of the basis elements of ${\cal Q}(R)$ on
basis elements of $p$-th order homogeneous component ${\cal M}^p(R)$
\be
(Q_1Q_{\overline 2}\dots Q_{\overline k})\triangleright
(M_{\overline{k+1}}M_{\overline {k+2}}\dots M_{\overline {k+p}}) =
\xi^k(M_{\underline{k+1}}M_{\underline {k+2}}\dots M_{\underline {k+p}}),
\quad \forall k,p\geqslant 1.
\label{Q-pol-ac}
\ee
Recall, that in the above formula  the copies of matrices are defined via
the R-matrix $R$:
$$
M_{\overline k} = R_{k-1}M_{\overline{k-1}}R^{-1}_{k-1},\quad
M_{\underline k} = R^{-1}_{k-1}M_{\underline{k-1}}R_{k-1}.
$$
\end{proposition}

Now, assuming $R$ to be a deformation of the flip $P$,  we discuss the restriction of
the adjoint BD algebra to some quotients of the RE algebra ${\cal M}(R)$ which can
be interpreted as a quantization of the coordinate algebra of $GL(m)$ orbits in $gl^*(m)$
\cite{GS1,GS3}. Such a quantum (braided) orbit is defined by an ideal $J_{\{c\}}$, generated
by elements
$$
{\Tr\str{-1.3}}_{\!R}(M^k) - c_k,\quad 1\leqslant k\leqslant m.
$$
(In order to get analogs of {\it regular} orbits we have to impose some restrictions on the
parameters $c_i$, see \cite{GS1}.)

The system of relations (\ref{N-M}) allows us to conclude, that the elements
${\Tr\str{-1.3}}_{\!R}(M^k)$ and ${\Tr\str{-1.3}}_{\!R}(N^k)$ are central in the adjoint BD
algebra. This is a consequence of the following property of the trace
$$
{\Tr\str{-1.3}}_{\!R^{(2)}}(R_1^{\pm 1}X_1R_1^{\mp 1}) =
{\Tr\str{-1.3}}_{\!R}(X)\,I_1,
$$
valid for an arbitrary $N\times N$ matrix $X$. Therefore, the quantum orbits are preserved
by the action of the RE algebra ${\cal Q}(R)$. Having restricted the adjoint BD algebra
${\cal B}_{\rm ad}({\cal Q},{\cal M})$ on the orbit ${\cal M}(R)/J_{\{c\}}$ we get the nontrivial
relations on the differential operators. They appear as the corresponding fixation of another
set of central elements --- ${\Tr\str{-1.3}}_{\!R}(N^k) = {\Tr\str{-1.3}}_{\!R} ((M^{-1}Q)^k)$.
\begin{definition}
A restriction of the adjoint BD algebra (\ref{Q-M}) on a quantum orbit
${\cal M}(R)/J_{\{c\}}$ is the quotient of ${\cal B}_{\rm ad}({\cal Q},{\cal M})$
over the ideal generated by the relations
\be
{\Tr\str{-1.3}}_{\!R}(M^k) = c_k,\qquad
{\Tr\str{-1.3}}_{\!R}((M^{-1}Q)^k) = {\Tr\str{-1.3}}_{\!R}((M^{-1}Q\,\triangleright)^k)
1_\sscr{$\,\cal B$}|_{J_{\{c\}}},\quad 1\leqslant k\leqslant m,
\label{restr}
\ee
where in the last relation we assume that traces of $M$ should be specified
to corresponding constants $c_i$ after calculating the action of $Q$.
\end{definition}
\begin{remark}
{\rm
The restriction $1\leqslant k\leqslant m$ in the above definition
is due to the fact that for a $GL(m)$-type R-matrix the quantum matrices $M$ and $M^{-1}Q$
satisfy the Cayley-Hamilton identity. The order of the Cayley-Hamilton polynomial is $m$,
so all traces ${\Tr\str{-1.3}}_{\!R}(M^p)$ and ${\Tr\str{-1.3}}_{\!R}((M^{-1}Q)^p)$ with
$p> m$ can be expressed in terms of the first $m$ traces.

Note also, that the restriction on central elements
${\Tr\str{-1.3}}_{\!R}((M^{-1}Q)^k)$ given in (\ref{restr}) is compatible
with the operator action of $Q$, presented in (\ref{Q-pol-ac}). Namely, one can show,
that ${\Tr\str{-1.3}}_{\!R}((M^{-1}Q\,\triangleright)^k)$ is a scalar operator on
any homogeneous component ${\cal M}^p(R)$. For example,
$$
{\Tr\str{-1.3}}_{\!R}(M^{-1}Q\,\triangleright) M_{1}M_{\overline {2}}\dots M_{\overline {p}} =
\xi\,{\Tr\str{-1.3}}_{\!R}(M^{-1})M_{1}M_{\overline {2}}\dots M_{\overline {p}}.
$$

At the classical level the corresponding restrictions have rather simple form $\Tr(M^kK) = 0$. These
relations mean that the $gl^*(m)$-module generated by infinitesimal vector fields arising from the
$GL(m)$ action on $gl^*(m)$ is a quotient of a free $gl^*(m)$-module. (However,  its restriction to a
generic orbit is projective.) The more complex formulae (\ref{restr}) are due to the ''exponentiated
character`` of quantum differential operators $Q$.}
\end{remark}

Consider a simple example, corresponding to a $GL(2)$-type R-matrix. Besides, we take the
generating matrices $M$ and $L$ to be of $2\times 2$ size.

The Cayley-Hamilton identity (\ref{CH-id}) for the matrix $M$ reads
$$
M^2 - qa_1(M)M +q^2a_2(M)I = 0,
$$
where
$$
a_1(M) = {\Tr\str{-1.3}}_{\!R}(M), \quad 2_qa_2(M) = q({\Tr\str{-1.3}}_{\!R}(M))^2
- {\Tr\str{-1.3}}_{\!R}(M^2)
$$
in accordance with (\ref{q-newt}).

Let us consider the "braided orbit" ${\cal O}(r)$ defined by the following values of the parameters
$c_k$, $k=1,\,2$:
\be
{\Tr\str{-1.3}}_{\!R}(M) = 0,\quad {\Tr\str{-1.3}}_{\!R}(M^2) = - \frac{2_q}{q^2}\,r^2,
\ee
where $r$ is a nonzero {\it real} number. Then the Cayley-Hamilton identity gives us the inverse
matrix $M^{-1}$ in the form
$$
M^{-1} = c\, M, \quad c = -r^{-2}.
$$
According to the definition (\ref{restr}), we can calculate the restriction conditions
for differential operators. The first condition is as follows
\be
{\Tr\str{-1.3}}_{\!R}(M^{-1}Q) = {\Tr\str{-1.3}}_{\!R}(M^{-1}Q\,\triangleright)
{1_\sscr{$\,\cal B$}}\raisebox{-1.5mm}{$|_{{\cal O}(r)}$}  =
\xi{\Tr\str{-1.3}}_{\!R}(M^{-1})\raisebox{-1.5mm}{$|_{{\cal O}(r)}$} =
\xi c\,{\Tr\str{-1.3}}_{\!R}(M)\raisebox{-1.5mm}{$|_{{\cal O}(r)}$} = 0,
\label{trMQ}
\ee
where we used (\ref{Q-unit}) for the action of $Q$ operator and the above explicit
form of $M^{-1}$ on the orbit ${\cal O}(r)$. Passing to the shifted set of generators
$Q_i^{\,j} = 1_\sscr{$\,\cal B$}\delta_i^j - (q-q^{-1})K_i^{\,j}$ we rewrite restriction
(\ref{trMQ}) in the form
\be
{\Tr\str{-1.3}}_{\!R}(MK) = 0.
\label{trMK}
\ee

At the classical limit $q\rightarrow 1$ the entries of the matrix $K$ become
generators of the Lie algebra $sl(2)$. We pass to the compact form of this
algebra by introducing new generators in the matrices $M$ and $K$. Namely,
we put
$$
M= \left(\begin{array}{cc}
ix_3&ix_1-x_2\\
ix_1+x_2&-ix_3
\end{array}\right).
$$
Also we have
\be
x_1^2+x_2^2+x_3^2 = r^2
\label{sphere}
\ee
as a consequence of the above equation  ${\Tr\str{-1.3}}_{\!R}(M^2) = -2_qq^{-2} r^2$.

Since at the limit $q\rightarrow 1$ $\Tr (K)$ is a central element, we can also take the matrix
$K$ to be traceless and parameterize its matrix elements as follows
$$
K=  \left(\begin{array}{cc}
iX_3&iX_1-X_2\\
iX_1+X_2&-iX_3
\end{array}\right).
$$
As follows from (\ref{Q-M}) the operators $X_i$ satisfy the $su(2)$ commutation relations:
$$
[X_i,X_j] = \varepsilon_{ijk}X_k,
$$
$\varepsilon_{ijk}$ being the completely antisymmetric tensor.

Taking into account the relations among $x_i$ and $X_i$
$$
[x_i,X_j] = - \varepsilon_{ijk}x_k
$$
(which can also be extracted as a classical limit of the third equation in (\ref{Q-M}))
we interpret generators $X_i$ as adjoint vector fields $X_i = \varepsilon_{ijk}x_j\partial_k$,
tangent to the sphere of radius $r$ described by the condition (\ref{sphere}). At the classical
limit the condition (\ref{trMK}) leads to the equality
$$
x_1X_1+x_2X_2+x_3X_3 = 0.
$$
This is the well-known identity which is satisfied by the tangent vector fields $X_i$.

We emphasize that  all higher order relations  following from (\ref{restr}) are satisfied automatically at
the classical limit. Indeed, calculating the second restriction for ${\Tr\str{-1.3}}_{\!R}((M^{-1}Q)^2)$
at the condition $\xi =1$ we get:
$$
{\Tr\str{-1.3}}_{\!R}(M^{-1}KM^{-1}+M^{-2}K - (q-q^{-1})M^{-1}KM^{-1}K)
 = {\Tr\str{-1.3}}_{\!R}(M^{-2}){\Tr\str{-1.3}}_{\!R}(I) -
({\Tr\str{-1.3}}_{\!R}(M^{-1}))^2.
$$
On the orbit ${\cal O}(r)$ we have ${\Tr\str{-1.3}}_{\!R}(M^{-1})=c\,{\Tr\str{-1.3}}_{\!R}(M)=0$,
$M^{-2} = c\,I$. So, at the classical limit we find:
$$
\Tr(MKM) = 2\,\Tr (M^2),
$$
or, using our parameterization for the matrices $M$ and $K$:
$$
\sum_{i,j,k}\varepsilon_{ijk}x_iX_jx_k = -2\sum_{i} x_i^2.
$$
The above relation indeed turns into an identity with the
choice $X_i= \varepsilon_{ijk}x_j\partial_k$.

\end{document}